\documentclass[11pt, reqno]{amsart}

\usepackage[english]{babel}

\usepackage{preamble}

\usepackage[a4paper,top=2cm,bottom=2cm,left=2.5cm,right=2.5cm,marginparwidth=2cm]{geometry}

\usepackage{amssymb, amsmath}
\usepackage{siunitx}
\PassOptionsToPackage{hyphens}{url}\usepackage{hyperref}
\usepackage{cleveref}
\usepackage[utf8]{inputenc}
\usepackage{csquotes}
\usepackage{booktabs}
\usepackage{longtable}
\usepackage{adjustbox}
\usepackage{array}
\usepackage{url}
\usepackage{xcolor} 
\usepackage{breqn}
\usepackage{tikz}
\usepackage{todonotes}

\usepackage[nobysame]{amsrefs}

\DeclareRobustCommand{\tfrac}{\textstyle\frac}

\title{Uniqueness for the Homogeneous Landau-Coulomb Equation in $L^{3/2}$}

\author{Maria Pia Gualdani}
\address{Department of Mathematics, The University of Texas at Austin, 2515 Speedway, Austin TX, 78712}
\email{gualdani@math.utexas.edu}

\author{Weiran Sun}
\address{Department of Mathematics, Simon Fraser University, 8888 University Dr., Burnaby, BC, Canada V5A 1S6}
\email{weirans@sfu.ca}

\date{}

\begin{document}

\begin{abstract}
We prove the uniqueness of $H$-solutions to the homogeneous Landau-Coulomb equation satisfying $\vint{v}^{k_0} f \in C([0, T]; L^{3/2}(\R^3))$ 
and $\vint{v}^{-3/2}\nabla_v ((\vint{v}^{k_0} f)^{3/4}) \in L^2((0, T) \times \R^3)$ 
for any $k_0 \geq 5$. In particular, this shows that the solutions constructed in~\cite{GGL25} are unique. The present work thus completes the global well-posedness theory in the critical space $L^{3/2}(\R^3)$. Our proof is part of a broader effort to use the $\CalM$-operator technique developed in~\cite{AGS2025, AMSY2020} to establish the  uniqueness of rough solutions to nonlinear kinetic equations. When applied to the space-homogeneous case, the $\CalM$-operator can be taken simply as a Bessel potential operator. 
\end{abstract}

\maketitle
\tableofcontents
\section{Introduction}


We consider the space-homogeneous Landau-Coulomb equation 
\begin{align}\label{eq:Landau-Coulomb}
  \partial_t f  
&= \nabla_v \vpran{\cdot {A[f] \nabla_v f - f \nabla_v a[f]}},
\qquad
  f|_{t=0} = f^{in}(v), 
\end{align}
with $A[f], \; a[f]$ being defined as
\begin{align} \label{def:A-a}
A[f] (t,v)& := \frac{1}{8\pi} \int_{\mathbb{R}^3} \frac{P(v-z) }{|v-z|} f(t,z)\dz,
\qquad 
 P(z) := I - \frac{z \otimes z}{|z|^2},  \nn
\\
a[f] (t,v) &:= Tr(A[f] ) = \frac{1}{4\pi} \int_{\mathbb{R}^3} \frac{1 }{|v-z|} f(t,z)\dz.
\end{align}
The goal of this paper is to establish the uniqueness of smooth solutions generated from rough initial data and the uniqueness of $H$-solutions satisfying additional regularity properties. 
Using the standard notation $\vint{v} := \sqrt{1 + |v|^2}$ for any $v \in \R^3$, our main result is summarized as follows.
\begin{main}
Let $k_0 \geq 5$, $\vint{v}^{k_0} f^{in} \in L^{3/2}(\R^3)$, and  $\|f^{in}\|_{L^1(\mathbb{R}^3)} =1 $. Let $T_0 > 0$ be arbitrary. If $f$ is a nonnegative $H$-solution to~\eqref{eq:Landau-Coulomb} that additionally satisfies
\begin{align} \label{cond:soln}
  \vint{v}^{k_0} f \in C([0, T_0]; L^{3/2}(\R^3)), 
\qquad
  \vint{v}^{-\frac{3}{2} + \frac{3k_0}{4}} \nabla_v (f^{3/4}) \in L^2((0, T_0) \times \R^3). 
\end{align}
Then it is unique.
\end{main}

\smallskip
Recall the definition of the $H$-solution (see Corollary 1.1 in \cite{D14}): for any test function $\phi \in C_c^2([0, T] \times \R^3)$, it holds that 
\begin{align*} 
& -\int_{\R^3} f^{in} \phi(0, v) \dv
- \int_0^T \int_{\R^3}
   f \del_t \phi \dv \dt
\\
& \hspace{3cm}
  = \int_0^T \int_{\R^3}
   A[f] f : \nabla^2 \phi \dv \dt 
  + \int_0^T \int_{\R^3}
     \nabla a[f] f \cdot \nabla \phi
     \dv \dt.
\end{align*}



\smallskip

Although the existence of smooth solutions to \eqref{eq:Landau-Coulomb} for initial data $f^{in} \in L^1(\mathbb{R}^3)$  (with suitable moment conditions) is well understood (see \cite{GS24, GGL25, JI24, DGGL24} and references therein), the question of uniqueness remains largely open. In particular, beyond the major unresolved problem of uniqueness for $H$-solutions, even the uniqueness of smooth solutions with initial data in $L^p$ with $1\le p \le \frac{3}{2}$, such as those constructed in \cite{GGL25, JI24, DGGL24}, has not yet been established. 

The Main Theorem above provides two contributions:  First, it shows that smooth solutions with initial data in $L^{3/2}$ are unique. The existence of such solutions was proved in \cite{GGL25}. Our result thus completes the global well-posedness analysis of the Landau equation with initial data in $L^p(\mathbb{R}^3)$ with $p \ge \frac{3}{2}$. 
Second, it yields a uniqueness result for $H$-solutions: if an $H$-solution additionally satisfies uniform bounds
in $L^{3/2}$ then it is unique. 

The conditional uniqueness of $H$-solutions is reminiscent of the work in  ~\cite{Fournier10}. In \cite{Fournier10}, uniqueness of $H$-solutions is obtained under a conditional uniform bound in $L^1(0, T; L^\infty(\R^3))$.  Many works have built their uniqueness theory  based on such $L^1(0, T; L^\infty(\R^3))$ bound, and therefore, by estimating the blow-up rate of $\norm{f(t)}_{L^\infty}$ at $t=0$ and showing that it is integrable (see, for example, \cite{GL2024,HW2024, SW2020}). In this regard, it has been shown in ~\cite{GGL25} that if $f^{in} \in L^p_{k_0}(\R^3)$ and $p > 3/2$, then the corresponding solution $f$ belongs to  $L^1(0, T; L^\infty(\R^3))$, and consequently, it is unique. Recently, an alternative proof of the uniqueness theorem in ~\cite{Fournier10}  has been provided in \cite{DGL25} using a novel estimate on the contractivity of the 2-Wasserstein distance. 

However, as shown in ~\cite{GGL25, JI24, DGGL24}, if the initial data is not smooth enough, one can often only expect $\norm{f(t)}_{L^\infty} \sim \frac{1}{t}$.  In this case, the uniqueness requires a novel approach and has been an open problem until now. 

Weak-strong uniqueness is another concept that is widely used to study kinetic equations (see, for example, \cite{Come2025, HST2020, HST2025, HW2024}). The basic idea is to show that if there exists a  solution with sufficient regularity, then any weak solutions with the same initial data will agree with the stronger solution. In this direction, the most recent result for the homogeneous Landau is in~\cite{Come2025}, where it is shown that the relative entropy between an $H$-solution $f$ and a regular solution $g$ is controlled by the initial relative entropy and the norms of $g$. The regularity of $g$ is stronger than $L^{3/2}$. 

The two routes mentioned above differ from our framework. There is one recent work~\cite{HJJ2025} that is closely related to ours, in which  negative Sobolev spaces and energy methods are used for the homogeneous Landau-Coulomb equation. The energy estimates are carried out using a Littlewood-Paley type decomposition and pseudo-differential analysis. The main uniqueness result in~\cite{HJJ2025} is established for initial data in a logarithmically modified Sobolev space $H^{-1/2, r}$. Note that in $\R^3$, the space $H^{-1/2}$ is precisely the $H^s$-space that $L^{3/2}$ embeds into. It is, however, unclear to us whether this embedding still holds for the modified space $H^{-1/2, r}$. In contrast, we directly work with $L^{3/2}$-spaces and our analysis is performed in the physical space.

\smallskip

\Ni {\bf Methodology.} 
The main idea in this work is to study uniqueness in suitable negative Sobolev spaces via direct energy estimates. To do so, we have been developing an $\CalM$-operator framework to tackle the issue of rough solutions.  We initiated our program in a preliminary work on spatial inhomogeneous kinetic equations in ~\cite{AGS2025}, where a model system and a modified Landau-Coulomb equation are shown to have unique bounded solutions. In~\cite{AGS2025},  we choose a  particular $\CalM$ with  symbol
\begin{align*}
  \CalM \sim 
\frac{1}{1 + \delta \int_t^T \vint{\xi + (t - \tau) \eta}^2 \dtau},
\end{align*}
where $(\xi, \eta)$ are the Fourier variables of $(v, x)$ respectively, and $\delta$ is a small parameter appropriately chosen. Such an operator was explicitly designed to have a manageable commutator with the transport operator $\del_t + v \cdot \nabla_x$.  The $\CalM$-operator defined above was first introduced in~\cite{AMSY2020} to show the regularization from $H^{-s}$ to $L^2$ of the linearized non-cutoff Boltzmann equation, where $s$ is the singularity strength of the Boltzmann operator.

The current work is part of this program and is an application of the $\CalM$-operator method to the space-homogeneous Landau-Coulomb equation. 
For space-homogeneous equations, the $\CalM$-operator above simplifies significantly.   We simply choose $\CalM$ to be the Bessel potential operator:
\begin{align*}
  \CalM = (I - \Delta_v)^{-1},
\end{align*}
where $I$ is the identity operator. The advantage is that such $\CalM$ has an integral representation in terms of the Bessel potential:
\begin{align} \label{integral_representation}
  (\CalM h)(v)
= \frac{1}{4\pi} \int_{\R^3}
   \frac{e^{-|v-v_1|}}{|v-v_1|}
   h(v_1) \dv_1. 
\end{align}
Although $\CalM$ has a simple symbol in the Fourier space, we choose to work in the physical space using its integral representation (\ref{integral_representation}). The main idea is to derive a closed $L^2$-energy estimate of $\CalM w$, where $w$ is the difference of two solutions generated from the same initial data (weighted with $\langle v \rangle^{2}$). Working with $\CalM w$ instead of $w$ allows one to consider less regular functions $w$. 
Unlike in ~\cite{AGS2025}, the particular choice of the Bessel potential avoids the need of pseudo-differential calculus to treat commutators. Instead, we use the following obvious yet convenient identity 
\begin{align} \label{reform:M-h}
  h = (I - \Delta_v) \CalM h,
\end{align}
and then work with the commutators generated by $I - \Delta_v$. This makes our computation explicit and elementary. 

Another important application of (\ref{integral_representation}) and ~\eqref{reform:M-h} is that it enables us to rewrite the nonlocal operators $A[w]$ and $a[w]$ in the Landau equation precisely into combinations of $A[\CalM w]$, $A[\Delta_v \CalM w]$, $a[\CalM w]$, and $a[\Delta_v \CalM w]$. Consequently, we obtain a closed-form equation for $\CalM w$. If one only relies on upper bounds and inequalities on $A[w]$ or $a[w]$, one might end up with terms of the form $\CalM |w|$ which cannot be controlled by $|\CalM w|$. We are currently investigating whether the techniques developed here for the homogeneous case can be used to improve the results in~\cite{AGS2025} for the inhomogeneous Landau equation. We also believe this $\CalM$-operator strategy could apply to the non-cutoff Boltzmann equations in both space-homogeneous and inhomogeneous cases. 

A final remark is in order. Our uniqueness result for $H$-solutions concerns the energy space associated with $C([0, T]; L^{3/2})$ and  the bound on $\nabla f^{3/4}$ as in the energy inequality. It is natural to ask whether uniqueness holds for $H$-solutions under the sole additional assumption that $f \in C([0, T]; L^{3/2})$. We show at the end of this paper that, in the {\it a priori} sense, the condition $f \in C([0, T]; L^{3/2})$ does indeed imply the desired bound on $\nabla f^{3/4}$. However, it is unclear whether $H$-solutions satisfy this {\it a priori} estimate. This may not be merely a technical limitation. In fact, as a comparison, in the case of incompressible Navier-Stokes, it is known that there are ``wild" weak solutions in $C(0, T; L^2)$ which are not Leray solutions \cite{BV2019}, even though in an {\it a priori} sense, the bound in $C(0, T; L^2)$ would yield a bound on the $L^2(0, T; H^1)$ norm. 

\medskip

{\bf{Notation:}} The following notations are adopted throughout this work. 
We denote $L^p(\R^n)$ and $W^{k,p}(\R^n)$ as the usual Lebesgue and Sobolev spaces. They are often abbreviated as $L^p$ and $W^{k,p}$ when there is no confusion. The weighted space $L^p_k(\R^3)$, or simply $L^p_k$, consists of functions satisfying
\begin{align*}
 h \in L^p_k(\R^3)
\Longleftrightarrow
 \vint{v}^k h \in L^p(\R^3).
\end{align*}
We reserve $A[\cdot]$ and $a[\cdot]$ as the nonlocal operators defined in~\eqref{def:A-a}. The operators $\nabla_v$, $\nabla_v \cdot$ and $\Delta_v$ always act on the velocity variable $v$, and the subscript $v$ is often omitted for brevity. The Hessian matrix of a function $h$ is denoted either as $Hess (h)$ or $\nabla^2 h$. The constant $C$ may change from line to line. The equivalence $\alpha \sim \beta$ is
\begin{align*}
  \alpha \sim \beta
\,\,\Longleftrightarrow \,\,
  c_1 \alpha \leq \beta \leq c_2 \alpha,
\end{align*}
for some $c_1, c_2 > 0$. Moreover, for any $\delta >0$, $\eta_\delta = \frac{1}{\delta^3}\eta_0 (x/ \delta)$ is the usual mollifier with support $B(0, \delta)$, $0 \leq \eta_0 \leq 1$, and $\int_{\R^3} \eta_0 (v) \dv = 1$.

\section{ Technical Lemmas}
Here we collect some technical lemmas that will be often used in later sections. 

\begin{lem} \label{ineq:basic} (Basic Inequalities)
Suppose $f, g$ are sufficiently regular such that each term in the inequalities is well-defined.

\Ni (a) Young's inequality:
\begin{align*}
  \norm{f \ast g}_{L^r}
\leq  
  \norm{f}_{L^p}
  \norm{g}_{L^q}, 
\qquad  
  1 + \frac{1}{r}
= \frac{1}{p} + \frac{1}{q}, 
\quad 
 1 \leq p, q \leq \infty.
\end{align*}

\Ni (b) $(L^p, L^q)$-interpolation:
\begin{align*}
 \norm{f}_{L^r}
\leq 
 \norm{f}_{L^p}^\theta
 \norm{f}_{L^q}^{1-\theta},
\qquad
  \frac{1}{r} 
= \frac{1}{p} \theta
  + \frac{1}{q} (1 - \theta),
\quad
  p < r < q,
\quad
  0 < \theta < 1.
\end{align*}

\Ni (c)  Hardy-Littlewood-Sobolev (HLS):
denote
\begin{align*}
 I_\lambda f 
= \int_{\R^n}\frac{f(y)}{|x - y|^\lambda} \dy, 
\qquad
  0 < \lambda < n. 
\end{align*}
Then
\begin{align*}
  \norm{I_\lambda f}_{L^q(\R^n)}
\leq 
  C \norm{f}_{L^p(\R^n)}, 
\qquad
  1 + \frac{1}{q}
= \frac{1}{p} + \frac{\lambda}{n}, 
\quad
1 < p < q < \infty. 
\end{align*}

\Ni (d) Sobolev embedding: if $1/p > k/n$, then 
\begin{align*}
 \dot{W}^{k,p}(\R^n) \hookrightarrow L^q(\R^n),
\qquad
\frac{1}{q} = \frac{1}{p} - \frac{k}{n}.
\end{align*}
If $1/p < k/n$, then
\begin{align*}
  W^{k,p}(\R^n)
\hookrightarrow
  C^{\ell, r}(\R^n),
\qquad
  \ell + \alpha = k - \frac{n}{p},
\qquad
  \alpha \in (0, 1), \,\, \ell = k - \left[\frac{n}{p}\right] - 1.
\end{align*}
 
\Ni (e) Bounds of the Bessel potential operator: 
\begin{align*}
  \norm{\CalM f}_{L^q(\R^n)}
&\leq  
  C \norm{f}_{L^p(\R^n)}, 
\qquad
  \frac{1}{p} - \frac{2}{n}
  \leq \frac{1}{q} \leq \frac{1}{p}, 
\quad
  p < \frac{n}{2},
\\
  \norm{\CalM^\beta f}_{L^p(\R^n)}
&\leq  
  C_p \norm{f}_{L^p(\R^n)}, 
\qquad \quad
  1 < p < \infty,
\qquad  
  \beta \geq 0,
\\
  \norm{\nabla_v^\alpha \CalM f}_{L^p(\R^n)}
&\leq  
  C_p \norm{f}_{L^p(\R^n)}, 
\qquad \quad
  1 < p < \infty,
\qquad
  |\alpha| \,\, = 1, 2.
\end{align*}

\Ni (f) Equivalent bounds:
\begin{align} \label{lem-ineq:f-1}
      C_2 \norm{\CalM \vpran{\vint{v}^{\alpha} f}}_{L^p(\R^3)}
\leq
  &\norm{\vint{v}^{\alpha} \CalM f}_{L^p(\R^3)}
\leq  
  C_1 \norm{\CalM \vpran{\vint{v}^{\alpha} f}}_{L^p(\R^3)},
\end{align}
for any $1 < p < \infty$ and $\alpha \in \R$. The constants $C_1, C_2$ only depend on $\alpha, p$. Moreover, it holds that
\begin{align} \label{lem-ineq:f-2}
  \norm{\vint{v}^\beta \nabla^{\alpha}\CalM f}_{L^p(\R^3)}
\leq
  C \norm{\vint{v}^\beta f}_{L^p(\R^3)},
\qquad
  |\alpha| \, = 1, 2,
\end{align}
for any $\beta \in \R $ and $1 < p < \infty$. 
\end{lem}
Bounds in $(a)-(d)$ are classical. Bounds in $(e)$ and $(f)$ are also classical in harmonic analysis. For completeness, we provide proofs of $(e)$ and $(f)$ in Appendix A. 

\smallskip
We will frequently use the bounds on $A [\cdot]$ and $a[\cdot]$ summarized in the lemma below.

\begin{lem} \label{lem:nabla-A-a}
({\bf A}) Let $F\ge 0$ be such that 
\begin{align*}
    c_1 \le \int F \;dv &\le c_2, \\
    \int F |v|^2 \;dv \le c_3, &\quad \int F \ln F \;dv \le c_4. 
\end{align*}
There exists $c_0 > 0$ depending on $c_1, c_2, c_3$ and $c_4$ such that
\begin{align} \label{cond:A-f}
  \xi \cdot A[F] \cdot \xi
\geq 
  c_0 \vint{v}^{-3} |\xi|^2,
\qquad
  \forall \, \xi \in \R^3.
\end{align}

\Ni ({\bf B}) Suppose $h$ is sufficiently smooth such that every term in the following inequalities is well-defined. Then, 

(a) for each $v \in \R^3$ and any $k \in \R$,
\begin{align*}
  \abs{\vint{v}^2 \nabla_v A [h]}
  + \abs{\vint{v}^2 \nabla_v a [h]}
& \leq
  C_k \vint{v}^{-1-k} \int_{|v-v_1| \leq \frac{1}{2} \vint{v}} 
 \frac{|\vint{v_1}^{3+k} h (v_1)|}{|v-v_1|^2} \dv_1
  + C \norm{h}_{L^1(\R^3)}. 
\end{align*}

(b) If, in addition, $\vint{v}^3 h  \in L^{\infty}(\R^3)$, then 
\begin{align*}
 \norm{\vint{\cdot}^2 \nabla_v A [h]}_{L^\infty(\R^3)}
 + \norm{\vint{\cdot}^2 \nabla_v a [h]}_{L^\infty(\R^3)}
\leq 
 C \norm{\vint{v}^3 h}_{L^\infty(\R^3)} + C \norm{h}_{L^1(\R^3)}.
\end{align*}

(c) for each $v \in \R^3$ and any $k \in \R$,
\begin{align*}
  \abs{A [h]} + \abs{a [h]}
& \leq
  C_k \vint{v}^{-k} \int_{|v-v_1| \leq \frac{1}{2} \vint{v}} 
 \frac{|\vint{v_1}^{1+k} h (v_1)|}{|v-v_1|^2} \dv_1
  + C \norm{h}_{L^1(\R^3)}. 
\end{align*}
\end{lem}

\begin{proof}
({\bf A}). The proof of the now-standard inequality (\ref{cond:A-f}) can be found in \cite{DV00} or in \cite{S17}.

\smallskip 

\Ni ({\bf B}). (a) By the definition of $A$, 
\begin{align*}
 \abs{\vint{v}^2 \nabla_v A [h]}
& \leq
 C\vint{v}^2 \int_{|v-v_1| \leq \frac{1}{2} \vint{v}} 
 \frac{|h(v_1)|}{|v-v_1|^2} \dv_1
 + C \vint{v}^2 \int_{|v-v_1| \geq 
 \frac{1}{2} \vint{v}}
 \frac{|h(v_1)|}{|v-v_1|^2} \dv_1
\\
&\leq
 C \vint{v}^{-1-k} \int_{|v-v_1| \leq \frac{1}{2} \vint{v}} 
 \frac{\vint{v_1}^{3+k}}{|v-v_1|^2} 
 \abs{h (v_1)} \dv_1
  + C \norm{h}_{L^1(\R^3)},
\end{align*}
where we have used that on the domain where $|v-v_1| \leq \frac{1}{2} \vint{v}$, it holds that 
\begin{align} \label{bound:equiv}
  c_2 \vint{v} \leq \vint{v_1} \leq c_1 \vint{v},
\end{align}
with $c_1, c_2$ being generic constants. The bound on $\vint{v}^2 \nabla a[h]$ holds similarly since $\nabla a$ and $\nabla A$ have similar upper bounds. 

\smallskip

\Ni (b) The proof follows immediately from the inequality in part (a) with $k=0$, since 
\begin{align*}
  \vint{v}^{-1} \int_{|v-v_1| \leq \frac{1}{2} \vint{v}} 
 \frac{1}{|v-v_1|^2}
 \abs{\vint{v_1}^3 h (v_1)} \dv_1
&\leq
  C \norm{\vint{\cdot}^3 h}_{L^\infty(\R^3)}.  
\end{align*}

\Ni (c) The proof is similar to part (a). By the definition of $A$,
\begin{align*}
  |A[h]|
&\leq
   C \int_{|v-v_1| \leq \frac{1}{2} \vint{v}} 
 \frac{|h(v_1)|}{|v-v_1|} \dv_1
 + C \int_{|v-v_1| \geq 
 \frac{1}{2} \vint{v}}
 \frac{|h(v_1)|}{|v-v_1|} \dv_1,
\end{align*}
where for $|v-v_1| \leq \frac{1}{2} \vint{v}$, we have $\frac{1}{|v-v_1|} \leq \frac{\vint{v}}{2}\frac{1}{|v-v_1|^2}$ and $\vint{v} \sim \vint{v_1}$. Therefore,
\begin{align*}
  |A[h]|
&\leq
   C \int_{|v-v_1| \leq \frac{1}{2} \vint{v}} 
 \frac{\vint{v}}{|v-v_1|^2} 
 \abs{h (v_1)} \dv_1
  + C \norm{h}_{L^1(\R^3)}
\\
& \leq
  C \vint{v}^{-k}
  \int_{|v-v_1| \leq \frac{1}{2} \vint{v}} 
 \frac{\vint{v_1}^{1+k}}{|v-v_1|^2} 
 \abs{h (v_1)} \dv_1
  + C \norm{h}_{L^1(\R^3)},
\end{align*}
where we have applied the equivalence in~\eqref{bound:equiv} over the domain $|v-v_1| \leq \frac{1}{2} \vint{v}$.
\end{proof}

In later proofs, we often need to separate $f$ into its smooth and non-smooth parts. The following lemma shows the uniform smallness of the non-smooth part of $f$. 
\begin{lem} \label{lem:unif-converg}
Suppose $f \in C([0, T_0]; L^{3/2}_{k_0}(\R^3))$. Let 
$f_\delta = f \ast \eta_\delta$. 
Then, for any $\epsilon > 0$, there exist $\delta_\ast, t_\ast > 0$ small enough such that for any $0 < \delta < \delta_\ast$,
\begin{align*}
  \Sup_{t \in [0, t_\ast]} \norm{\vint{v}^{k_0} \vpran{f(t) - f_\delta(t)}}_{L^{3/2}(\R^3)}
< \epsilon.
\end{align*}
This also implies that  
\begin{align} \label{bound:unif-f-delta}
  \Sup_{t \in [0, t_\ast]}\norm{\vint{v}^{k_0} f_\delta}_{L^{3/2}}
\leq  
 C_1 < \infty, 
\end{align}
where $C_1$ is independent of $\delta$.
\end{lem}
\begin{proof}
First, by the continuity in time, there exists $t_\ast$ small enough such that 
\begin{align*}
  \Sup_{t \in [0, t_\ast]} \norm{\vint{v}^{k_0} f(t) - \vint{v}^{k_0} f^{in}}_{L^{3/2}(\R^3)}
< \epsilon/4.
\end{align*}
Young's inequality gives that for each $t \leq t_\ast$,  
\begin{align*}
  \norm{(\vint{\cdot}^{k_0} f) \ast \eta_\delta - (\vint{v}^{k_0} f^{in}) \ast \eta_\delta}_{L^{3/2}(\R^3)}
\leq
 \norm{\vint{v}^{k_0} f(t) - \vint{v}^{k_0} f^{in}}_{L^{3/2}(\R^3)}
< \epsilon/4. 
\end{align*}
Choose $\delta_\ast < 1$ small enough such that for any $0 < \delta < \delta_\ast < 1$, 
\begin{align*}
  \norm{\vint{v}^{k_0} f^{in} - (\vint{v}^{k_0} f^{in}) \ast \eta_\delta}_{L^{3/2}(\R^3)}
\leq 
  \epsilon/4. 
\end{align*}
Furthermore, 
\begin{align*}
  \abs{\vint{v}^{k_0} f_\delta
  - (\vint{v}^{k_0} f) \ast \eta_\delta}
&\leq 
  \int_{\R^3} \abs{\vint{v}^{k_0} - \vint{v_1}^{k_0}} |f(v_1)|
  \eta_0 \vpran{\frac{v-v_1}{\delta}} \dv_1
\\
& \leq
 C \int_{|v-v_1| \leq \delta} |v_1 - v| \abs{\vint{v_1}^{k_0-1} f(v_1)}
 \eta_0 \vpran{\frac{v-v_1}{\delta}} \dv_1
\\
& \leq
 C \delta \vpran{\vint{\cdot}^{k_0-1} |f|} \ast \eta_\delta,
\end{align*}
where we have used that $|v_1 - v| \leq \delta < 1$, and consequently, $|v_1 - v| \leq 2 \vint{v_1}$. By Young's inequality, there exists $C > 0$ independent of $t, \delta$ such that for each $t \in [0, t_\ast]$, 
\begin{align*}
  \norm{\vint{v}^{k_0} f_\delta
  - (\vint{v}^{k_0} f_\delta) \ast \eta_\delta}_{L^{3/2}(\R^3)}
& \leq 
 C \delta \norm{\vpran{\vint{\cdot}^{k_0-1} |f|} \ast \eta_\delta}_{L^{3/2}(\R^3)}
\\
& \leq
 C \delta \norm{\vint{\cdot}^{k_0-1} f}_{L^{3/2}(\R^3)}
\leq 
 C \delta
< \epsilon/4,
\end{align*}
by taking $\delta < \frac{\epsilon}{4C}$. 
Therefore, by applying the triangle inequality, we obtain that 
\begin{align*}
& \quad \,
  \Sup_{t \in [0, t_\ast]} \norm{\vint{v}^{k_0} f(t) - \vint{v}^{k_0} f_\delta(t) }_{L^{3/2}(\R^3)}
\\
&\leq
 \Sup_{t \in [0, t_\ast]} \norm{\vint{v}^{k_0} f(t) - \vint{v}^{k_0} f^{in}}_{L^{3/2}(\R^3)}
\\
& \quad \,
 + \Sup_{t \in [0, t_\ast]} \norm{\vpran{\vint{v}^{k_0} f} \ast \eta_\delta(t) - \vpran{\vint{v}^{k_0} f^{in}} \ast \eta_\delta}_{L^{3/2}(\R^3)}
\\
& \quad \,
 + \Sup_{t \in [0, t_\ast]} \norm{\vpran{\vint{v}^{k_0} f^{in}} \ast \eta_\delta - \vint{v}^{k_0} f^{in}}_{L^{3/2}(\R^3)}
\\
& \quad \,
 + \Sup_{t \in [0, t_\ast]}
   \norm{\vint{v}^{k_0} f_\delta
  - (\vint{v}^{k_0} f) \ast \eta_\delta}_{L^{3/2}(\R^3)}
< \epsilon. \qedhere
\end{align*}
\end{proof}

\begin{rmk} \label{rmk:unif-Lp}
Note that the space $C([0, T_0]; L^{3/2}_{k_0}(\R^3))$ in Lemma~\ref{lem:unif-converg} can be replaced by any $C([0, T_0]; L^{p}_{k_0}(\R^3))$ with $1 < p < \infty$. 
\end{rmk}
From now on, we will restrict the time interval to $[0, T] \subseteq [0, t_\ast] \subseteq [0, T_0]$ and $0 < \delta < \delta_\ast$.

\section{{\it A Priori} Estimates} \label{Sec:a-priori}
In this section, we perform the {\it a priori} estimates for the Landau equation.
More specifically, suppose $f, g$ are two solutions to~\eqref{eq:Landau-Coulomb} with the same initial data $f^{in}$ and $f, g$ satisfy~\eqref{cond:soln}. 
Denote
\begin{align*}
  w_0 := f - g,
\qquad
  w := \vint{v}^2 (f - g).
\end{align*}
Then $w$ satisfies the linear equation
\begin{align} \label{eq:w}
 \del_t w
&= \nabla_v \cdot (A [f] \nabla_v w) - \nabla_v \cdot (w \nabla a[f])
+ \vint{v}^2 
\nabla_v \cdot (A [w_0] \nabla_v g) 
\\
& \quad \, 
- \vint{v}^2 \nabla_v \cdot (g \nabla a[w_0]) \nn
  + R_1 + R_2, 
\end{align}
with $w^{in} = 0$ and 
\begin{align} \label{def:R-12-orig}
 & R_1
 = \vint{v}^2 \nabla_v \cdot (A [f] \nabla_v w_0) - \nabla_v \cdot (A [f] \nabla_v w),
\\
 & R_2
 = \nabla_v \cdot (w \nabla a[f])
   - \vint{v}^2 \nabla_v \cdot (w_0 \nabla_v a[f]).
\end{align}
We can simplify $R_1$ and $R_2$ as
\begin{align} \label{def:R-12}
  R_1
  = -2v \cdot A [f] \nabla_v w_0
    - 2 \nabla_v \cdot (A[f] \, v \, w_0),
\qquad
   R_2
  = 2 \, v \, w_0 \cdot \nabla_v a[f].
\end{align}
Multiplying ~\eqref{eq:w} by the test function $\phi = \CalM^2 w$ and integrating in $(t, v)$ 
we get
\begin{align} \label{decomp:CalM-w}
  \norm{\CalM w(T)}^2_{L^2(\R^3)}
&= \int_0^T \int_{\R^3} 
  \CalM w \, \nabla_v \CalM \cdot
  (A[f] \nabla_v w) \dv\dt  \nn
\\
& \quad \,
 - \int_0^T \int_{\R^3} 
  \CalM^2 w \, \nabla_v  \cdot
  (w \nabla_v a[f]) \dv\dt  \nn
\\
& \quad \, 
  + \int_0^T \int_{\R^3} 
  \CalM^2 w \, \vint{v}^2 \nabla_v \cdot
  (A[w_0] \nabla_v g) \dv\dt  \nn 
\\
& \quad \,
 - \int_0^T \int_{\R^3} 
  \CalM^2 w \, \vint{v}^2 \nabla_v \cdot
  (g \nabla_v a[w_0]) \dv\dt \nn
\\
& \quad \,
 + \int_0^T \int_{\R^3} 
  \CalM^2 w \, R_1 \dv\dt \nn
 + \int_0^T \int_{\R^3} 
  \CalM^2 w \, R_2 \dv\dt
\\
&=: \text{Int}(I_1) + \text{Int}(I_2) + \text{Int}(I_3) + \text{Int}(I_4) + \text{Int}(I_5) + \text{Int}(I_6).
\end{align}
In the following subsections, we estimate each $I_k$ separately. As mentioned in the introduction, to convert estimates on $w, w_0$ into bounds on $\CalM w, \CalM w_0$, we will often use the identities  
\begin{align} \label{reform:w-Mw}
   w = (I - \Delta_v) \CalM w,
\qquad
   w_0 = (I - \Delta_v) \CalM w_0.
\end{align}

\subsection{Estimate of $\mathrm{Int}(I_1)$} Denote
\begin{align*}
  I_1 
:= \int_{\R^3} 
  \CalM w \, \nabla_v  \cdot
 \CalM( (A[f] \nabla_v w)) \dv, 
\qquad 
  \text{Int}(I_1)
= \int_0^T I_1 \dt.
\end{align*}
Through integration by parts and~\eqref{reform:w-Mw}, we decompose $I_1$ as follows: 
\begin{align*}
  I_1 
:= & \int_{\R^3} 
  \CalM w \, \nabla_v  \cdot
 \CalM( (A[f] \nabla_v w)) \dv \\
=& - \int_{\R^3} 
  \nabla_v \CalM w \cdot
 \CalM( (A[f] \nabla_v w)) \dv \\
 =& - \int_{\R^3} 
  \nabla_v \CalM w \cdot
 \CalM( A[f] \CalM^{-1} \nabla_v \CalM w) \dv ,
\end{align*}
using the fact that $ \nabla_v w = \nabla_v ( \CalM^{-1} \CalM w )=\CalM^{-1}\nabla_v   \CalM w  $ . 
From $I_1$, we will get a purely coercive term 
$$D : =- \int_{\R^3} \langle 
  \nabla_v \CalM w,  
   A[f] \nabla \CalM w  \rangle \dv$$ 
   and lower order terms. For that, we first have to commute $\CalM^{-1}$ with $A[f]$: recall that $\CalM^{-1} = I - \Delta_v$ and 
   $$
   h_1 \; \CalM^{-1} h_2 = \CalM^{-1}(h_1 h_2) + h_2 \Delta_v h_1  + 2 \nabla_v h_1 \cdot \nabla_v h_2,
   $$
   which yields 
   \begin{align*}
   \CalM( A[f] \CalM^{-1} \nabla_v \CalM w) = \;& A[f]  \nabla_v \CalM w + \CalM (  \Delta_v A[f] \nabla_v \CalM w) 
   +  2 \CalM (\partial_{v_i} A[f]  \nabla_v \partial_{v_i} \CalM w).
   \end{align*}
Therefore,
\begin{align*}
I_1 =& - \int_{\R^3} \langle 
  \nabla_v \CalM w,  
   A[f] \nabla \CalM w  \rangle \dv \\
  &  - \int_{\R^3} \langle 
  \nabla_v \CalM w, \CalM (  \Delta_v A[f] \nabla_v \CalM w) \rangle \dv \\
 &  - 2 \int_{\R^3} \langle 
  \nabla_v \CalM w,\CalM (\partial_{v_i} A[f]  \nabla_v \partial_{v_i} \CalM w)  \rangle \dv \\
  =&  - \int_{\R^3} \langle 
  \nabla_v \CalM w,  \;
   A[f] \nabla \CalM w  \rangle \dv \\
  &  - \int_{\R^3} \langle 
  \nabla_v \CalM^2  w,  \;  \Delta_v A[f] \nabla_v \CalM w  \rangle \dv \\
 &  - 2 \int_{\R^3} \langle 
  \nabla_v \CalM^2  w,\;  \partial_{v_i} A[f]  \nabla_v \partial_{v_i} \CalM w \rangle \dv.
        \end{align*}
        Integration by parts in the last integral yields 
        \begin{align*}
I_1 =&  - \int_{\R^3} \langle 
  \nabla_v \CalM w,  
   A[f] \nabla \CalM w  \rangle \dv \\
  &  + 2 \int_{\R^3} \langle 
  \partial_{v_i} \nabla_v \CalM^2  w,\;  \partial_{v_i} A[f]  \nabla_v  \CalM w \rangle \dv  \\
  &  + \int_{\R^3} \langle 
  \nabla_v \CalM^2  w,  \;  \Delta_v A[f] \nabla_v \CalM w  \rangle \dv \\
   =: &  \; D + I_{1,1} + I_{1,2}.
  \end{align*}

The dissipation term $D$ gives
\begin{align} \label{est:D}
  D
\leq 
  - c_0 \int_0^T \norm{\vint{v}^{-3/2} \nabla_v \CalM w}^2_{L^2(\R^3)} \dt, 
\end{align}
where $c_0$ depends on the mass, second moment, and entropy of the function $f$.

\smallskip
\Ni \underline{\it Estimates of $\mathrm{Int} (I_{1,2})$} 
In $\mathrm{Int} (I_{1,2})$, we first remove the gradient from $ \nabla_v \CalM w$. Integrating by parts, we rewrite $I_{1,2}$ as
\begin{align} 
  I_{1,2}
= & - \int_{\R^3}
    \nabla_v \cdot \vpran{ \partial_{v_i v_i}A[ f] \nabla_v \CalM^2 w } \CalM w \dv \nn \\
     =& - \int_{\R^3} \nabla_v \partial_{v_i v_i} a[ f] \cdot \nabla_v \CalM^2 w \CalM w \dv  - \int_{\R^3} \CalM w  \; Tr(  \partial_{v_i v_i} A[ f]  Hess( \CalM^2 w)) \dv \label{decomp:I-12-interm}
\\ 
    =& \int_{\R^3} \nabla_v f  \cdot \nabla_v \CalM^2 w \CalM w \dv  - \int_{\R^3} \CalM w  \; Tr(  \partial_{v_i v_i} A[ f]  Hess( \CalM^2 w)) \dv \nn,
    \end{align}
where we have used
\begin{align*} 
\nabla_v \cdot A[h] = \nabla_v a [h], \qquad
  - h
=  \Delta_v  a[h].
\end{align*} 
In the second term of $I_{1,2}$, we now remove one derivative from $A[f]$, yielding: 
\begin{align*}
   \int_{\R^3} \CalM w    \; Tr(  \partial_{v_i v_i} A[ f]  \; Hess( \CalM^2 w)) \dv 
    = & - \int_{\R^3} \partial_{v_i} \CalM w  \; Tr(  \partial_{v_i} A[ f]\;   Hess( \CalM^2 w)) \dv \\
    & - \int_{\R^3}  \CalM w  \; Tr(  \partial_{v_i} A[ f]  \; Hess( \partial_{v_i} \CalM^2 w)) \dv .
\end{align*}
In summary, $I_{1,2}$ is the sum of three terms of the form
\begin{align} 
  I_{1,2} =   & - \int_{\R^3}  \CalM w  \; Tr(  \partial_{v_i} A[ f]  \; Hess( \partial_{v_i} \CalM^2 w)) \dv \nn\\
 &  - \int_{\R^3} \partial_{v_i} \CalM w  \; Tr(  \partial_{v_i} A[ f]\;   Hess( \CalM^2 w)) \dv  \nn\\
  & + \int_{\R^3} (\nabla_v f  \cdot \nabla_v \CalM^2 w)\; \CalM w \dv
  =: \;  I_{1,2}^{(1)} + I_{1,2}^{(2)}
   + I_{1,2}^{(3)}.
   \label{decomp:I-12}
  \end{align}


We start with $I_{1,2}^{(1)}$. 
\begin{lem} \label{lem:I-121}
Let $\epsilon_0 > 0$ be an arbitrarily small number. Then there exists $C_{\epsilon_0} > 0$ such that 
\begin{align*}
 \mathrm{Int}(I_{1,2}^{(1)})
\leq 
 \epsilon_0 \int_0^T \norm{\vint{v}^{-3/2} \nabla_v \CalM w}_{L^2}^2 \dt
+ C_{\epsilon_0} \int_0^T \norm{\CalM w}_{L^2}^2 \dt. 
\end{align*}
\end{lem}
\begin{proof}
To deal with the term $ \partial_{v_i}  A[f]$ we separate $f$ into the smooth and less smooth parts, as highlighted in Lemma~\ref{lem:unif-converg}. For that,  fix any $\epsilon > 0$ small. Let $\delta_\ast$ be defined in Lemma~\ref{lem:unif-converg} and denote $\tilde f = f - f_\delta$ with $0 < \delta < \delta_\ast$. Then we have
\begin{align*}
  f = f_\delta + \tilde f, 
\qquad
  \Sup_{t \in [0, t_\ast]} \norm{\tilde f (t)}_{L^{3/2}(\R^3)}
< \epsilon. 
\end{align*}
Then 
\begin{align*}
  \mathrm{Int} (I_{1,2}^{(1)})
=&   - \int_{\R^3}  \CalM w  \; Tr(  A[\partial_{v_i} f_\delta]  \; Hess( \partial_{v_i} \CalM^2 w)) \dv 
\\
& - \int_{\R^3}  \CalM w  \; Tr(  \partial_{v_i}A[ \tilde f]  \; Hess( \partial_{v_i} \CalM^2 w)) \dv
\\
=&:    \mathrm{Int} (I_{1,2}^{(1, 1)})
   + \mathrm{Int} (I_{1,2}^{(1, 2)}). 
\end{align*}
The integral with the smooth $f_\delta$ is bounded as follows:
\begin{align*}
  \abs{\mathrm{Int} (I_{1,2}^{(1, 1)})}
&\leq 
 C \int_0^T \int_{\R^3}
     \abs{\vint{v}^{-3/2} Hess( \partial_{v_i} \CalM^2 w))} \,  
    \abs{\vint{v}^{3/2} \partial_{v_i} A[f_\delta]} \, \abs{\CalM w} \dv \dt  \nn \\
  &\leq  C \vpran{\norm{\vint{v}^3 f_\delta}_{L^\infty} + 1}\int_0^T \int_{\R^3}\abs{\vint{v}^{-3/2} Hess( \partial_{v_i} \CalM^2 w)} \, \abs{\CalM w} \dv \dt.  \nn \\
  \end{align*}
Commuting $\vint{v}^{-3/2}$ with $ Hess$ in the first $L^2$-norm, we get 
\begin{align*}
    \|{\vint{v}^{-3/2} Hess( \partial_{v_i} \CalM^2 w))}\|_{L^2} \; \le & \; \|Hess (\vint{v}^{-3/2}  \partial_{v_i} \CalM^2 w)\|_{L^2}\\
    & + \| \vint{v}^{-3/2} \partial_{v_i} \CalM^2 w\|_{L^2} +  \|  \partial_{v_i}\partial_{v_j} \CalM^2 w\|_{L^2} \\
    =& \; \|\Delta_v (\vint{v}^{-3/2}  \CalM (\partial_{v_i} \CalM w))\|_{L^2} \\& + \| \vint{v}^{-3/2} \partial_{v_i} \CalM^2 w\|_{L^2} + \|  \Delta \CalM^2 w\|_{L^2}.
\end{align*}
We have the following bounds: since $\Delta = -\CalM^{-1}+ I $, we get
\begin{align}
\|  \Delta \CalM^2 w\|_{L^2} = \| \CalM^2 w\|_{L^2} + \|   \CalM w\|_{L^2} \le 2 \|   \CalM w\|_{L^2},\label{1}
\end{align}
using the second inequality in Lemma  \ref{ineq:basic} with $\beta =1$ and $p=2$. 
Moreover 
\begin{align}
     \| \vint{v}^{-3/2} \partial_{v_i} \CalM^2 w\|_{L^2} = \| \vint{v}^{-3/2} \CalM(\partial_{v_i} \CalM w)\|_{L^2} \le \| \vint{v}^{-3/2} \partial_{v_i} \CalM w\|_{L^2},\label{2}
\end{align}
thanks to Lemma  \ref{ineq:basic}$(f)$ and the second inequality of Lemma  \ref{ineq:basic}$(e)$ (with $\beta =1$). 
To estimate the term with the Laplacian, we first observe that 
\begin{align*}
    \Delta_v (\vint{v}^{-3/2}  (\partial_{v_i} \CalM^2 w)) = & \; \Delta_v (\vint{v}^{-3/2})  \partial_{v_i} \CalM^2 w \\
    &+ \vint{v}^{-3/2}  \Delta_v(\partial_{v_i} \CalM^2 w) \\
    &+ 2 \nabla_v \vint{v}^{-3/2} \cdot \nabla_v (\partial_{v_i} \CalM^2 w) \\
     = & \; \Delta_v (\vint{v}^{-3/2}) \CalM  (\partial_{v_i} \CalM w) \\
    &- \vint{v}^{-3/2}  \partial_{v_i} \CalM w + \vint{v}^{-3/2} \CalM (\partial_{v_i} \CalM w)  \\
    &+{{2 \nabla_v \vint{v}^{-3/2} \cdot \nabla_v (\partial_{v_i} \CalM^2 w)}}.
\end{align*}
Each term can be bounded as 
\begin{align}
 \|\Delta_v (\vint{v}^{-3/2}) \CalM  (\partial_{v_i} \CalM w) \|_{L^2}\le \|   \vint{v}^{-3/2} \partial_{v_i} \CalM w\|_{L^2}, \nn\\
 \|\vint{v}^{-3/2} \CalM (\partial_{v_i} \CalM w) _{L^2}\le \|   \vint{v}^{-3/2} \partial_{v_i} \CalM w\|_{L^2}, \label{3}\\
 \|{2 \nabla_v \vint{v}^{-3/2} \cdot \nabla_v (\partial_{v_i} \CalM^2 w)}\|_{L^2} \le \| Hess(\CalM^2 w)\|_{L^2} \le \|\CalM w\|_{L^2}, \nn
 \end{align}
where in the last inequality we apply the third estimate of Lemma~\ref{ineq:basic} part (e) with $\alpha =2$. 
Summarizing, from (\ref{1}), (\ref{2}), and (\ref{3}) we obtain 
\begin{align*}
    \|{\vint{v}^{-3/2} Hess( \partial_{v_i} \CalM^2 w))}\|_{L^2} \; \lesssim \|\CalM w\|_{L^2} + \|   \vint{v}^{-3/2} \partial_{v_i} \CalM w\|_{L^2}.
    \end{align*}
The bound above and Young's inequality applied  to $I_{1,2}^{(1, 1)}$ yield 
\begin{align} \label{est:I-1211}
  \abs{\mathrm{Int} (I_{1,2}^{(1, 1)})}
&\leq 
 C \vpran{\norm{\vint{v}^3 f_\delta}_{L^\infty} + 1}\int_0^T \int_{\R^3}\abs{\vint{v}^{-3/2} Hess( \partial_{v_i} \CalM^2 w)} \, \abs{\CalM w} \dv \dt  \nn  \\
 &\le \epsilon_0 \int_0^T  \norm{   \vint{v}^{-3/2} \partial_{v_i} \CalM w}^2_{L^2}\;dt + \frac{C\vpran{\norm{\vint{v}^3 f_\delta}_{L^\infty} + 1}^2}{\epsilon_0} \int_0^T  \norm{\CalM w}^2_{L^2}\;dt .
  \end{align}

Moreover, we have that 
\begin{align*}
  \norm{\vint{v}^3 f_\delta}_{L^\infty}
\leq 
  \frac{1}{\delta^3}
  \norm{\vint{\cdot}^3 f}_{L^{3/2}(\R^3)}
  \vpran{\int_{\R^3} \eta_0^3 \vpran{\frac{z}{\delta}} \dz}^{1/3}
\leq  
  \frac{C}{\delta^2}
  \norm{\vint{\cdot}^3 f}_{L^{3/2}(\R^3)}.
\end{align*}


Next, we estimate $\mathrm{Int}(I_{1,2}^{(1,2)})$. Using Lemma~\ref{lem:nabla-A-a} part (a) with $k = 1/2$, we have
\begin{align} \label{boudn:I-12-12-interm}
  \abs{\mathrm{Int} \vpran{I_{1,2}^{(1,2)}}}
&\leq
 C \int_0^T \int_{\R^3}
 \abs{\vint{v}^{-3/2} \nabla^2_v \CalM \vpran{\nabla_v\CalM w}} \,  
    \abs{\vint{v}^{3/2} \nabla_v A[\tilde f]} \, \abs{\CalM w} \dv \dt \nn
\\
& \hspace{-0.5cm}\leq
  C \int_0^T \int_{\R^3}
  \abs{\vint{v}^{-3/2} \nabla^2_v \CalM \vpran{\nabla_v\CalM w}} \,  
    \vpran{\int_{\R^3} 
 \frac{|\vint{v_1}^{7/2} \tilde f (v_1)|}{|v-v_1|^2} \dv_1} \, \abs{\vint{v}^{-3/2}\CalM w} \dv \dt \nn
\\
& \quad \,
  + C \int_0^T \norm{\tilde f}_{L^1(\R^3)} \int_{\R^3}
  \abs{\vint{v}^{-3/2} \nabla^2_v \CalM \vpran{\nabla_v\CalM w}} \,  
     \abs{\CalM w} \dv \dt 
\\
& \hspace{-0.5cm} \leq
  C \int_0^T 
    \norm{\vint{v}^{-3/2} \nabla^2_v \CalM \vpran{\nabla_v\CalM w}}_{L^2(\R^3)}
    \norm{\int_{\R^3} 
 \frac{|\vint{v_1}^{7/2} \tilde f (v_1)|}{|v-v_1|^2} \dv_1}_{L^3(\R^3)} \nn
\\
& \hspace{8.5cm} \times
 \norm{\vint{v}^{-3/2}\CalM w}_{L^6(\R^3)} \dt \nn
\\
& \quad \, 
  + C \norm{f}_{L^\infty_t(L^1(\R^3))}
  \int_0^T \int_{\R^3}
  \abs{\vint{v}^{-3/2} \nabla_v^2 \CalM \vpran{\nabla_v\CalM w}} \abs{\CalM w} \dv \dt.  \nn
\end{align}
In the last term, we apply Young's inequality and then (\ref{lem-ineq:f-2})  with $\alpha = 2$ and $\beta = -\frac{3}{2}$ to  get 
\begin{align}\label{est11}
C \norm{f}_{L^1(\R^3)}
  &\int_0^T \int_{\R^3}
  \abs{\vint{v}^{-3/2} \nabla_v^2 \CalM \vpran{\nabla_v\CalM w}} \abs{\CalM w} \dv \dt \nn\\
  &\le \frac{\epsilon_0}{4}
  \int_0^T \norm{\vint{v}^{-3/2}\vpran{\nabla_v\CalM w}}_{L^2(\R^3)}^2 \dt + \frac{C}{\epsilon_0} \int_0^T \norm{\CalM w}_{L^2(\R^3)}^2 \dt, 
\end{align}
where $C$ depends on the $L^1$ norm of $f$.

For the other term, HLS inequality yields
\begin{align*}
  \norm{\int_{\R^3} 
 \frac{|\vint{v_1}^{7/2} \tilde f (v_1)|}{|v-v_1|^2} \dv_1}_{L^3(\R^3)}
\leq
  \norm{\vint{\cdot}^{7/2} \tilde f}_{L^{3/2}}. 
\end{align*}
Hence, by the Sobolev embedding of $\dot H^1(\R^3) \hookrightarrow L^6(\R^3)$ and (\ref{lem-ineq:f-2}) we get 
\begin{align}
C \int_0^T 
    \norm{\vint{v}^{-3/2} \nabla^2_v \CalM \vpran{\nabla_v\CalM w}}_{L^2(\R^3)}&
    \norm{\int_{\R^3} 
 \frac{|\vint{v_1}^{7/2} \tilde f (v_1)|}{|v-v_1|^2} \dv_1}_{L^3(\R^3)} \nn
\\
&  \times
 \norm{\vint{v}^{-3/2}\CalM w}_{L^6(\R^3)} \dt \nn \\
  \leq
  \vpran{C \sup_{[0, T]}  \norm{\vint{\cdot}^{7/2} \tilde f}_{L^{3/2}} }
  &\int_0^T \norm{\vint{v}^{-3/2}\vpran{\nabla_v\CalM w}}_{L^2(\R^3)}^2 \dt \label{est12}
\end{align}
Combining (\ref{est11}) and (\ref{est12}) we get
\begin{align*}
  \abs{\mathrm{Int} \vpran{I_{1,2}^{(1,2)}}}
&\leq
  \vpran{C \sup_{[0, T]}  \norm{\vint{\cdot}^{7/2} \tilde f}_{L^{3/2}} + \frac{\epsilon_0}{4}}
  \int_0^T \norm{\vint{v}^{-3/2}\vpran{\nabla_v\CalM w}}_{L^2(\R^3)}^2 \dt
\\
& \quad \quad  \,
 + \frac{C}{\epsilon_0} \int_0^T \norm{\CalM w}_{L^2(\R^3)}^2 \dt. 
\end{align*}
By Lemma~\ref{lem:unif-converg}, we can choose $\delta$ small enough such that
\begin{align*}
  C \sup_{[0, T]}  \norm{\vint{\cdot}^{7/2} \tilde f}_{L^{3/2}}
 = C \sup_{[0, T]}  \norm{\vint{\cdot}^{7/2} \vpran{f - f_\delta}}_{L^{3/2}}
\leq \epsilon_0/4,
\end{align*}
which gives
\begin{align} \label{est:I-1212}
  \abs{\mathrm{Int} \vpran{I_{1,2}^{(1,2)}}}
\leq
  \frac{\epsilon_0}{2}
  \int_0^T \norm{\vint{v}^{-3/2}\vpran{\nabla_v\CalM w}}_{L^2(\R^3)}^2 \dt
 + C \int_0^T \norm{\CalM w}_{L^2(\R^3)}^2 \dt. 
\end{align}
Combining~\eqref{est:I-1211} and~\eqref{est:I-1212}, we obtain the desired bound for $\mathrm{Int} (I^{(1)}_{1,2})$. 
\end{proof}

Next, we estimate $\mathrm{Int} (I^{(2)}_{1,2})$.

\begin{lem} \label{lem:I-122}
Let $\epsilon_0 > 0$ be an arbitrarily small number. There exists $C_{\epsilon_0} > 0$ such that 
\begin{align*}
 \mathrm{Int}(I_{1,2}^{(2)})
\leq 
 \epsilon_0 \int_0^T \norm{\vint{v}^{-3/2} \nabla_v \CalM w}_{L^2}^2 \dt
+ C_{\epsilon_0} \int_0^T \norm{\CalM w}_{L^2}^2 \dt. 
\end{align*}
\end{lem}

\begin{proof}
This term can be handled with similar techniques as in $\mathrm{Int} (I^{(1)}_{1,2})$.  The only difference is that the operator $\nabla^2 \CalM$ acts on $\nabla \CalM w$ instead of  $\CalM w$.  We have
\begin{align*}
  \abs{\mathrm{Int} (I_{1,2}^{(2)})}
&\leq 
 C \int_0^T \int_{\R^3}
     \abs{\vint{v}^{-3/2} \vpran{\nabla_v\CalM w}} \,  
    \abs{\vint{v}^{3/2} \nabla_v A[f]} \, \abs{\nabla^2_v \CalM  \CalM w} \dv \dt. 
\end{align*}

For the smooth part of $f$, we obtain an estimate similar to (\ref{est:I-1211}) after using (\ref{lem-ineq:f-2}) to estimate $\|\nabla^2_v \CalM  \CalM w\|_{L^2} $. Similar estimates for $\mathrm{Int} (I^{(1)}_{1,2})$ can be applied to bound the term with the non smooth part of $f$. The resulting bound is similar to (\ref{est:I-1212}).  The details are henceforth omitted to avoid repetition.
\end{proof}

The last term to be estimated in the decomposition \eqref{decomp:I-12} is  $\mathrm{Int} (I_{1,2}^{(3)})$. This is the content of the next lemma. 

\begin{lem} \label{lem:I-123}
Let $\epsilon_0 > 0$ be an arbitrarily small number. Then there exists $C_{\epsilon_0} > 0$ such that 
\begin{align*}
 \mathrm{Int}(I_{1,2}^{(3)})
\leq 
 \epsilon_0 \int_0^T \norm{\vint{v}^{-3/2} \nabla_v \CalM w}_{L^2}^2 \dt
 + \epsilon_0 \, \Sup_{[0, T]} \norm{\CalM w}_{L^2}^2
 + C_{\epsilon_0} \int_0^T \norm{\CalM w}_{L^2}^2 \dt. 
\end{align*}
\end{lem}

\begin{proof}
Recall 
$$
I_{1,2}^{(3)}:=  \int_{\R^3} (\nabla_v f  \cdot \nabla_v \CalM^2 w)\; \CalM w \dv.
 $$
We first write 
$$
\nabla_v f = \tfrac{4}{3}\; f^{\frac{1}{4}} \; \nabla_v f^{\frac{3}{4}}.
$$
We separate $\nabla_v f^{3/4} = \nabla_v (f^{3/4}) \ast \eta_\delta + h$, where
\begin{align} \label{def:h}
  h (t, x) = \nabla_v (f^{3/4}) - \nabla_ v (f^{3/4}) \ast \eta_\delta. 
\end{align}
Note that $f^{3/4} \in L^\infty(0, T; L^2(\R^3))$ and $\nabla_ v f^{3/4} \in L^2((0, T) \times \R^3)$. Together with Young's inequality, this implies that
\begin{align*}
  \nabla_v(f^{3/4}) \ast \eta_\delta
\in L^2((0, T) \times \R^3), 
\qquad
  h \in L^2((0, T) \times \R^3).
\end{align*}
First, we show that $h$ is small in the sense of~\eqref{limit:h}.

We claim that $
  \lim_{\delta \to 0}
  \int_0^T \int_{\R^3} |h|^2 \dv \dt
= 0. $ Let $\hat h$ be the Fourier transform in $v$, then
\begin{align*}
  \int_0^T \int_{\R^3} |h|^2 \dv \dt 
= \int_0^T \int_{\R^3}
  |1 - \hat\eta_0(\delta \xi)|^2  \, |\xi|^2 |\widehat {\vpran{f^{3/4}}} (t, \xi)|^2 \dxi \dt. 
\end{align*}
Since $\eta_\delta = \frac{1}{\delta^3} \eta_0(x/\delta)$ is the smooth mollifier with compact support, we have 
\begin{align*}
  \hat \eta_0 (0) = 1,
\qquad
  \hat \eta_0 \in C^\infty(\R^3),
\qquad
  \norm{\hat \eta_0}_{L^\infty(\R^3)} \leq \norm{\eta_0}_{L^1(\R^3)} = 1.
\end{align*}
This shows that for each $t, \xi$, 
\begin{align*}
  |1 - \hat\eta_0(\delta \xi)|^2  \, |\xi|^2 |\widehat {\vpran{f^{3/4}}} (t, \xi)|^2 \to 0 
\qquad
  \text{as $\delta \to 0$}, 
\end{align*}
and 
\begin{align*}
  |\hat h|^2
\, \leq \, 
  4 |\xi|^2 |\widehat {\vpran{f^{3/4}}}|^2 
\, \in 
  L^1((0, T) \times \R^3). 
\end{align*}
Therefore, the Lebesgue Dominated Convergence theorem applies, and we have 
\begin{align} \label{limit:h}
  \lim_{\delta \to 0}
  \int_0^T \int_{\R^3} |h|^2 \dv \dt
= 0. 
\end{align}

With this in mind, we rewrite $\mathrm{Int}(I_{1,2}^{(3)})$ as
\begin{align*}
  \mathrm{Int}(I_{1,2}^{(3)})
& = \tfrac{4}{3} \int_{\R^3}
    \vpran{\nabla_v \CalM^2 w \cdot 
    f^{\frac{1}{4}} \nabla_v f^{\frac{3}{4}}} \CalM w \dv
\\
& = - \tfrac{4}{3} \int_{\R^3}
    \vpran{\nabla_v \CalM^2 w \cdot 
    f^{\frac{1}{4}} \vpran{\nabla_v f^{\frac{3}{4}} \ast \eta_\delta}} \CalM w \dv
  - \tfrac{4}{3} \int_{\R^3}
    \vpran{\nabla_v \CalM^2 w \cdot 
    f^{\frac{1}{4}} h} \CalM w \dv
\\
& =: \mathrm{Int}(I_{1,2}^{(3,1)})
    + \mathrm{Int}(I_{1,2}^{(3,2)}), 
\end{align*}
where $h$ is defined in~\eqref{def:h}. 
We bound $\mathrm{Int}(I_{1,2}^{(3,1)})$ by
\begin{align} \label{est:I-1231}
  \abs{\mathrm{Int}(I_{1,2}^{(3,1)})}
&\leq 
  \int_0^T \norm{\nabla (\CalM^2 w)}_{L^6}
  \norm{f^{\frac{1}{4}}}_{L^6}
  \norm{\nabla_v f^{\frac{3}{4}} \ast \eta_\delta}_{L^2}
  \norm{\CalM w}_{L^6}\dt \nn
\\
&\leq 
  \frac{C}{\delta}\vpran{\Sup_{[0, T]} \norm{f}_{L^{\frac{3}{2}}}}
  \int_0^T \norm{\CalM w}_{L^2}
  \norm{\nabla \CalM w}_{L^2}\dt\nn
\\
& \leq
  \frac{\epsilon_0}{4}
  \int_0^T
  \norm{\nabla \CalM w}^2_{L^2}\dt
  + \frac{C}{\delta^2} 
  \vpran{\Sup_{[0, T]} \norm{f}_{L^{\frac{3}{2}}}}^2
  \int_0^T \norm{\CalM w}^2_{L^2} \dt, 
\end{align}
where we have applied the Sobolev inequality and Young's inequality:
\begin{align*}
  &\norm{\nabla (\CalM^2 w)}_{L^6(\R^3)}
\leq 
  C \norm{\CalM^2 w}_{L^2(\R^3)}\leq 
  C \norm{\CalM w}_{L^2(\R^3)},
\quad
 \norm{\CalM w}_{L^6(\R^3)}
\leq
  C_1 \norm{\nabla \CalM w}_{L^2(\R^3)}, 
\end{align*}
and
\begin{align*}
  \norm{\nabla_v f^{\frac{3}{4}} \ast \eta_\delta}_{L^2}
= \norm{f^{\frac{3}{4}} \ast \nabla \eta_\delta}_{L^2}
\leq 
  \frac{C}{\delta}
  \norm{f^{\frac{3}{4}}}_{L^2(\R^3)}
= \frac{C}{\delta}
  \norm{f}^{\frac{3}{4}}_{L^{\frac{3}{2}}(\R^3)}. 
\end{align*}
The second term $\mathrm{Int}(I_{1,2}^{(3,2)})$ satisfies
\begin{align*}
  \abs{\mathrm{Int}(I_{1,2}^{(3,2)})}
&\leq 
  C \int_0^T \norm{\nabla\CalM (\CalM w)}_{L^6}
  \norm{f^{\frac{1}{4}}}_{L^6}
  \norm{h}_{L^2}
  \norm{\CalM w}_{L^6}\dt  
\\
& \leq 
 C \vpran{\Sup_{[0, T]} \norm{f}_{L^{\frac{3}{2}}}}
 \int_0^T 
 \norm{\CalM w}_{L^2}
 \norm{h}_{L^2}
 \norm{\nabla \CalM w}_{L^2} \dt.
\end{align*}
Therefore, by the limit~\eqref{limit:h}, we choose $\delta$ sufficiently small such that
\begin{align} \label{est:I-1232}
& \quad \, 
  \abs{\mathrm{Int}(I_{1,2}^{(3,2)})} \nn
\\
& \leq
  C \vpran{\Sup_{[0, T]} \norm{f}_{L^{\frac{3}{2}}}}
 \vpran{\Sup_{[0, T]} \norm{\CalM w}_{L^2}} 
 \vpran{\int_0^T \norm{h}_{L^2}^2 \dt}^{\frac{1}{2}}
 \vpran{\int_0^T \norm{\nabla \CalM w}_{L^2}^2 \dt}^{\frac{1}{2}} \nn
\\
& \leq
  C \vpran{\int_0^T \norm{h}_{L^2}^2 \dt} \vpran{\Sup_{[0, T]} \norm{\CalM w}_{L^2}}^2
 +   C \vpran{\int_0^T \norm{h}_{L^2}^2 \dt} 
 \vpran{\int_0^T \norm{\nabla \CalM w}_{L^2}^2 \dt} \nn
\\
& \leq
 \frac{\epsilon_0}{2} 
 \Sup_{[0, T]} \norm{\CalM w}^2_{L^2}
 + \frac{\epsilon_0}{2}
 \int_0^T \norm{\nabla \CalM w}_{L^2}^2 \dt.
\end{align}
The desired bound for $\mathrm{Int}(I_{1,2}^{(3)})$ is obtained by combining~\eqref{est:I-1231} and~\eqref{est:I-1232}. 
\end{proof}

By adding the bounds in Lemmas~\ref{lem:I-121}, ~\ref{lem:I-122}, ~\ref{lem:I-123} and redefining $3 \epsilon_0$ as $\epsilon_0$, we obtain the bound for $\mathrm{Int}(I_{1,2})$ as follows. 
\begin{lem} \label{lem:I-12}
Let $\epsilon_0 > 0$ be an arbitrarily small number. Then there exists $C_{\epsilon_0} > 0$ such that 
\begin{align*}
 \abs{\mathrm{Int}(I_{1,2})}
\leq 
 \epsilon_0 \int_0^T \norm{\vint{v}^{-3/2} \nabla_v \CalM w}_{L^2}^2 \dt
 + \epsilon_0 \, \Sup_{[0, T]} \norm{\CalM w}_{L^2}^2
 + C_{\epsilon_0} \int_0^T \norm{\CalM w}_{L^2}^2 \dt. 
\end{align*}    
\end{lem}

With the estimates above, we can now show the bound for $\mathrm{Int}(I_1)$ defined in~\eqref{decomp:CalM-w}. 

\begin{prop} \label{prop:I-1}
Let $\epsilon_0 > 0$ be an arbitrarily small number. Then there exists $C_{\epsilon_0} > 0$ such that 
\begin{align*}
 \mathrm{Int}(I_{1})
\leq 
 {\color{blue}{(-{c_0}+\epsilon_0})} \int_0^T \norm{\vint{v}^{-3/2} \nabla_v \CalM w}_{L^2}^2 \dt
 + \epsilon_0 \, \Sup_{[0, T]} \norm{\CalM w}_{L^2}^2
 + C_{\epsilon_0} \int_0^T \norm{\CalM w}_{L^2}^2 \dt. 
\end{align*}   
\end{prop}
\begin{proof}
It remains to  bound $\mathrm{Int}(I_{1,1})$. Note that $I_{1,1}$ is similar to $\mathrm{Int}(I_{1,2}^{(2)})$, since
\begin{align*}
  \mathrm{Int}(I_{1,1}) 
= 2 \int_0^T \int_{\R^3}
   (\nabla_v^2 \CalM^2 w : \nabla A [f]) \cdot \nabla_v \CalM w \dv \dt
= -2 \, \mathrm{Int}(I_{1,2}^{(2)}).  
\end{align*}
Thus the bound for $\mathrm{Int}(I_1)$ follows from Lemma~\ref{lem:I-12} and the lower bound for $D$ in (\ref{est:D}). 
\end{proof}

\subsection{Estimates for $\mathrm{Int}(I_2)$} 
Recall that 
\begin{align*}
  \mathrm{Int}(I_{2})
= - \int_0^T \int_{\R^3} 
  \CalM^2 w \, \nabla_v \cdot
  (w \nabla_v a[f]) \dv\dt. 
\end{align*}


\begin{prop} \label{prop:I-2}
Let $\epsilon_0 > 0$ be an arbitrarily small number. Then there exists $C_{\epsilon_0} > 0$ such that 
\begin{align*}
 \abs{ \mathrm{Int}(I_{2})}
\leq 
 \epsilon_0 \int_0^T \norm{\vint{v}^{-3/2} \nabla_v \CalM w}_{L^2}^2 \dt
 + \epsilon_0 \, \Sup_{[0, T]} \norm{\CalM w}_{L^2}^2
 + C_{\epsilon_0} \int_0^T \norm{\CalM w}_{L^2}^2 \dt. 
\end{align*}  
\end{prop}

\begin{proof}

We substitute  $w$  with $ (I - \Delta) \CalM w$ and integrate by parts in $\mathrm{Int}(I_{2})$ to get
\begin{align*}
 \mathrm{Int}(I_{2})
&= \int_0^T \int_{\R^3}
  \nabla_v \CalM^2 w \cdot \nabla_v a[f] \,
  (I - \Delta_v) \CalM w \dv \dt
\\
& = \int_0^T \int_{\R^3}
  (I - \Delta_v)
  \vpran{\nabla_v \CalM^2 w \cdot \nabla_v a[f]} \CalM w \dv \dt
\\
& = \int_0^T \int_{\R^3}
  (I - \Delta_v) (\nabla_v \CalM^2 w)
  \cdot \nabla_v a[f] \, \CalM w \dv \dt
\\
& \quad \, 
  -2 \int_0^T \int_{\R^3}
  \nabla_v (\partial_{v_i}\CalM^2 w) \cdot 
  \nabla_v  (\partial_{v_i} a[f]) \, \CalM w \dv \dt
\\
& \quad \,
  + \int_0^T \int_{\R^3}
    \nabla_v \CalM^2 w
    \cdot (-\Delta_v) \nabla_v a[f] \, \CalM w \dv \dt
\\
&=: \mathrm{Int}(I_{2,1})
   + \mathrm{Int}(I_{2,2})
   + \mathrm{Int}(I_{2,3}).
\end{align*}
The term $\mathrm{Int}(I_{2,1})$ 
is of a similar structure as $$\mathrm{Int}(I_{1,2}^{(1)})= - \int_{\R^3}    \; Tr(  A[\partial_{v_i} f]  \; Hess( \partial_{v_i} \CalM^2 w)) \CalM w\dv,$$
with $Hess( \partial_{v_i} \CalM^2 w)$ replaced by $(I - \Delta_v) (\nabla_v \CalM^2 w)$  and $\partial_{v_i} A[ f]$ by $\nabla_v a[f]$. Note that the sup norm of $\partial_{v_i} A[ f]$ and of  $\nabla_v a[f]$ are equivalent, as well as the $L^2$ norm of $Hess( \partial_{v_i} \CalM^2 w)$ and of $\Delta_v(\nabla_v \CalM^2 w)$. Therefore, without repeating the details, we conclude that 
\begin{align*}
 |\mathrm{Int}(I_{2,1})|
\leq 
 \epsilon_0 \int_0^T \norm{\vint{v}^{-3/2} \nabla_v \CalM w}_{L^2}^2 \dt
+ C_{\epsilon_0} \int_0^T \norm{\CalM w}_{L^2}^2 \dt. 
\end{align*}


The term $\mathrm{Int}(I_{2,3})$ equals $\mathrm{Int}(I_{1,2}^{(3)})$ since $\Delta_v a[f] = -f$, which implies 
\begin{align*}
 |\mathrm{Int}(I_{2,3})|
\leq 
 \epsilon_0 \int_0^T \norm{\vint{v}^{-3/2} \nabla_v \CalM w}_{L^2}^2 \dt
 + \epsilon_0 \, \Sup_{[0, T]} \norm{\CalM w}_{L^2}^2
 + C_{\epsilon_0} \int_0^T \norm{\CalM w}_{L^2}^2 \dt. 
\end{align*}



It remains to estimate  $\mathrm{Int}(I_{2,2})$. Integration by parts yields 
\begin{align*}
  \mathrm{Int}(I_{2,2})
& = 2 \int_0^T \int_{\R^3}
  \partial_{v_i}  \vpran{\CalM w \nabla_v (\partial_{v_i}  \CalM^2 w)} \cdot
  \nabla_v a[f] \,  \dv \dt
\\
& \hspace{-0.5cm}
 = 2 \int_0^T \int_{\R^3}
  \vpran{\Delta_v \nabla_v \CalM^2 w} \cdot
  \nabla_v a[f] \, \CalM w \,  +
  \partial_{v_i} (\CalM w) \nabla_v (\partial_{v_i}  \CalM^2 w) \cdot
  \nabla_v a[f] \,  \dv \dt 
\\
& \hspace{-0.5cm}
  = 2 \int_0^T \int_{\R^3}
  \vpran{\nabla_v \CalM^2 w} \cdot
  \nabla_v a[f] \, \CalM w \,  \dv \dt 
  - 2 \int_0^T \int_{\R^3}
  \vpran{\nabla_v \CalM w} \cdot
  \nabla_v a[f] \, \CalM w \,  \dv \dt
\\
& \quad \,
  + 2 \int_0^T \int_{\R^3}
  \partial_{v_i} (\CalM w) \nabla_v (\partial_{v_i}  \CalM^2 w) \cdot
  \nabla_v a[f] \dv \dt.
\end{align*}
The first two terms  are similar to $\mathrm{Int}(I_{1,2}^{(1)})$ and the third term is similar to $\mathrm{Int}(I_{1,2}^{(2)})$. Therefore, all the terms in $\mathrm{Int}(I_2)$ are bounded similarly to the terms in $\mathrm{Int}(I_1)$ and the desired bound holds. 
\end{proof}

\begin{rmk}
Alternatively, the term $\mathrm{Int}(I_{2,1})$ can be bounded using the following procedure. We have
\begin{align*}
  \mathrm{Int}(I_{2,1})
&= \int_0^T \int_{\R^3}
  (\nabla_v \CalM w)
  \cdot \nabla_v a[f_\delta] \, \CalM w \dv \dt
  + \int_0^T \int_{\R^3}
  (\nabla_v \CalM w)
  \cdot \nabla_v a[\tilde f] \, \CalM w \dv \dt.
\end{align*}
First the first term we use Lemma~\ref{lem:nabla-A-a} part (b) and get: 
\begin{align*}
& \quad \,
 \abs{\int_0^T \int_{\R^3}
  (\nabla_v \CalM w)
  \cdot \nabla_v a[f_\delta] \, \CalM w \dv \dt}
\\
&\leq 
 C \Sup_{[0, T]} \norm{\vint{v}^{3/2}\nabla_v a[f_\delta]}_{L^\infty(\R^3)}
 \vpran{\int_0^T
  \norm{\vint{v}^{-3/2} \nabla_v \CalM w}^2_{L^2} \dt}^{1/2}
 \vpran{\int_0^T \norm{\CalM w}^2_{L^2} \dt}^{1/2}
\\
& \leq
C \Sup_{[0, T]}
 \norm{\vint{v}^3 f_\delta}_{L^\infty}
 \vpran{\int_0^T
  \norm{\vint{v}^{-3/2} \nabla_v \CalM w}^2_{L^2} \dt}^{1/2}
 \vpran{\int_0^T \norm{\CalM w}^2_{L^2} \dt}^{1/2}
\\
&\leq
 \frac{C}{\delta^2} \Sup_{[0, T]}\norm{f}_{L^{3/2}}
 \vpran{\int_0^T
  \norm{\vint{v}^{-3/2} \nabla_v \CalM w}^2_{L^2} \dt}^{1/2}
 \vpran{\int_0^T \norm{\CalM w}^2_{L^2} \dt}^{1/2} 
\\
& \leq
 \epsilon_0 \int_0^T
  \norm{\vint{v}^{-3/2} \nabla_v \CalM w}^2_{L^2} \dt
  + \frac{C}{\delta^4} 
  \int_0^T \norm{\CalM w}^2_{L^2} \dt. 
\end{align*}
Meanwhile, the second term in $\mathrm{Int}(I_{2,1})$ satisfies
\begin{align*}
& \quad \, 
  \int_0^T \int_{\R^3}
  (\nabla_v \CalM w)
  \cdot \nabla_v a[\tilde f] \, \CalM w \dv \dt
\\
&= \frac{1}{2} \int_0^T \int_{\R^3}
   \abs{\CalM w}^2 (-\Delta_v) a[\tilde f] \dv \dt
 = \frac{1}{2} \int_0^T \int_{\R^3}
   \abs{\CalM w}^2 \, \tilde f\dv \dt
\\
& \leq
 \int_0^T 
   \norm{\CalM w}^2_{L^6}
   \norm{\tilde f}_{L^{3/2}} \dt
 \leq
 \int_0^T 
   \norm{\nabla \CalM w}^2_{L^2}
   \norm{\tilde f}_{L^{3/2}} \dt
\\
& \leq 
 \Sup_{[0, T]}\norm{\tilde f}_{L^{3/2}}
 \int_0^T 
   \norm{\nabla \CalM w}^2_{L^2} \dt
\leq
  \epsilon_0 \int_0^T 
   \norm{\nabla \CalM w}^2_{L^2} \dt,
\end{align*}
by taking $\delta$ small enough. Adding these two parts again gives that 
\begin{align*}
  \abs{\mathrm{Int}(I_{2,1})}
\leq
  \epsilon_0 \int_0^T 
   \norm{\nabla \CalM w}^2_{L^2} \dt
  + C_{\epsilon_0}
  \int_0^T 
   \norm{\nabla \CalM w}^2_{L^2} \dt.
\end{align*}
\end{rmk}

\medskip

\subsection{Estimates of $\mathrm{Int}(I_3)$ and $\mathrm{Int}(I_4)$} Unlike in $I_1$ or $I_2$, the unweighted difference $f - g$ in $I_3$ and $I_4$  appears inside the nonlocal operators $A[\cdot]$ or $a[\cdot]$. The  estimates required  for these terms therefore slightly differ from those used so far. In fact, the need for weighted estimates arises precisely from these contributions: we will control $A[f-g]$ and $\nabla a[f-g]$ by $L^2$  norms of $\CalM w$, where $w$ was defined as $\langle v \rangle ^{2}(f-g)$. This is the content of the next proposition. 

\begin{prop} \label{prop:A-divA-w}
Let $ w_0 = f-g$. The terms $A[w_0]$ 
and $\nabla_v a[w_0]$ have the following expressions in terms of $\CalM w_0$ and $\nabla_v \CalM w_0$:
\begin{align} \label{reform:A-w}
  A[w_0]
= A[\CalM w_0]
    - A[\Delta \CalM w_0],
\end{align}
with the bounds
\begin{align} \label{bound:A-Mw0}
  \abs{A[\CalM w_0] (v)}
\leq
  C \vint{v} \norm{\CalM w}_{L^2(\R^3)},
\end{align}
and
\begin{align} \label{bound:A-Delta-Mw0}
  \norm{A[\Delta \CalM w_0]}_{L^6(\R^3)}
\leq
  C \norm{\vint{\cdot}^{-3/2}\nabla \CalM w}_{L^2(\R^3)}
  + C \norm{\CalM w}_{L^2(\R^3)},
\end{align}
where recall that $w = \vint{v}^2 w_0$.
Moreover, 
\begin{align}
 \nabla_v \cdot & A[w_0]
= \nabla_v a[w_0] 
= A_1
  + \nabla_v \CalM w_0. \label{reform:div-A-w}
\end{align}
where $A_1$ satisfies
\begin{align} \label{bound:A-1}
  \norm{A_1}_{L^6(\R^3)}
\leq
  C \norm{\CalM w_0}_{L^2(\R^3)}.
\end{align}
\end{prop}
\begin{proof}
The key step is again to write $w_0 = (I - \Delta) \CalM w_0$. One has
\begin{align*}
  \nabla_v \cdot A[w_0]
&= \nabla_v a[w_0]
= \frac{1}{4 \pi}
  \nabla_v \int_{\mathbb{R}^3} \frac{1}{|v-v_1|} (I - \Delta_{v_1}) \CalM w_0(v_1) \dv_1
\\
&= \frac{1}{4 \pi} \nabla_v 
  \int_{\R^3} (I - \Delta_{v_1}) \vpran{\frac{1}{|v-v_1|}} \CalM w_0(v_1) \dv_1
\\
& = \frac{1}{4 \pi} \nabla_v 
  \int_{\R^3} \frac{1}{|v-v_1|} \CalM w_0(v_1) \dv_1
  + \nabla_v 
  \int_{\R^3} \delta_0(v-v_1) \CalM w_0(v_1) \dv_1
\\
& = \frac{1}{4 \pi} \nabla_v 
  \int_{\R^3} \frac{1}{|v-v_1|} \CalM w_0(v_1) \dv_1
  + \nabla_v \CalM w_0. 
\end{align*}
Let 
\begin{align*}
  A_1
= \frac{1}{4 \pi} \nabla_v 
  \int_{\R^3} \frac{1}{|v-v_1|} \CalM w_0(v_1) \dv_1.
\end{align*}
Then by HLS, $A_1$ satisfies
\begin{align*}
  \norm{A_1}_{L^6(\R^3)}
\leq
  C \norm{\int_{\R^3} \frac{1}{|v-v_1|^2} \CalM w_0(v_1) \dv_1}_{L^6(\R^3)}
\leq
 C \norm{\CalM w_0}_{L^2(\R^3)}.
\end{align*}
Similarly, 
\begin{align*}
  A[w_0]
&= \frac{1}{8 \pi} 
  \int_{\R^3}
  \frac{P(v-v_1)}{|v-v_1|}
  (I - \Delta_{v_1}) \CalM w_0(v_1) \dv_1
 = A[\CalM w_0]
  - A[\Delta \CalM w_0].  
\end{align*}
Separating the integration domain, we have
\begin{align*}
  \abs{A[\CalM w_0] (v)}
&\leq
  C \int_{\R^3}
  \frac{1}{|v-v_1|}
  \abs{\CalM w_0(v_1)} \dv_1
\\
& \hspace{-1cm} \leq
  C \int_{|v-v_1| \leq \frac{1}{2} \vint{v}}
  \frac{1}{|v-v_1|}
  \abs{\CalM w_0(v_1)} \dv_1
  +  C \int_{|v-v_1| \geq \frac{1}{2} \vint{v}}
  \frac{1}{|v-v_1|}
  \abs{\CalM w_0(v_1)} \dv_1. 
\end{align*}
Applying Cauchy-Schwarz to the first term gives
\begin{align*}
  \int_{|v-v_1| \leq \frac{1}{2} \vint{v}}
  \frac{1}{|v-v_1|}
  \abs{\CalM w_0(v_1)} \dv_1
& \leq
  C \vint{v} \norm{\CalM w_0}_{L^2(\R^3)}
\leq
  C \vint{v} \norm{\CalM w}_{L^2(\R^3)},
\end{align*}
where the last inequality follows from Lemma~\ref{ineq:basic} part $(f)$. Moreover, it holds that
\begin{align*}
 \int_{|v-v_1| \geq \frac{1}{2} \vint{v}}
  \frac{1}{|v-v_1|}
  \abs{\CalM w_0(v_1)} \dv_1
&\leq
 C \vint{v}^{-1} \norm{\CalM w_0}_{L^1(\R^3)}
\\
& \hspace{-2cm} \leq
 C \vint{v}^{-1}
 \norm{\vint{v}^2\CalM w_0}_{L^2(\R^3)}
\leq
 C \vint{v}^{-1}
 \norm{\CalM w}_{L^2(\R^3)}. 
\end{align*}
Adding these two bounds gives the desired bound for $\abs{A[\CalM w_0] (v)}$. Note that the embedding of $\CalM w_0$ from $L^1$ to a weighted-$L^2$ is the only reason that we need to consider a weighted energy estimate. 

Finally, to bound $A[\Delta \CalM w_0]$, we use its definition and obtain
\begin{align*}
  \abs{A[\Delta \CalM w_0]}
&= \frac{1}{8 \pi}
  \abs{\int_{\R^3}
  \frac{P(v-v_1)}{|v-v_1|}
  \Delta_{v_1} \CalM w_0(v_1) \dv_1}
\\
& = \frac{1}{8 \pi}
  \abs{\int_{\R^3}
  \nabla_{v_1}\vpran{\frac{P(v-v_1)}{|v-v_1|}}
  \cdot \nabla_{v_1} \CalM w_0(v_1) \dv_1}
\\
& \leq
  C \int_{\R^3}
  \frac{1}{|v - v_1|^2}
  \abs{\nabla_{v_1} \CalM w_0(v_1)} \dv_1.
\end{align*}
By HLS, 
\begin{align*}
 \norm{A[\Delta \CalM w_0]}_{L^6(\R^3)}
\leq
 C \norm{\nabla \CalM w_0}_{L^2(\R^3)}
\leq
 C \norm{\vint{\cdot}^{-3/2} \nabla \CalM w}_{L^2(\R^3)}
 + C \norm{\CalM w}_{L^2(\R^3)},
\end{align*}
where the last inequality follows from Lemma~\ref{ineq:basic} part $(f)$.
\end{proof}

Now we show the bounds of $\mathrm{Int}(I_3)$ and $\mathrm{Int}(I_4)$.

\begin{prop} \label{prop:I-3-I-4}
Let $\epsilon_0 > 0$ be an arbitrarily small number. Then there exists $C_{\epsilon_0} > 0$ such that 
\begin{align*}
 \abs{ \mathrm{Int}(I_{3})}
\leq 
    \epsilon_0 \, \Sup_{[0, T]}
  \norm{\CalM w}^2_{L^2}
  + \epsilon_0 
  \int_0^T
  \norm{\vint{v}^{-3/2} \nabla_v \CalM w}^2_{L^2} \dt  \nn
  + C_{\epsilon_0} \int_0^T 
  \norm{\CalM w}^2_{L^2(\R^3)} \dt.
\end{align*}  
A similar bound holds for $\mathrm{Int}(I_{4})$.
\end{prop}
\begin{proof}
First, we introduce some notations. Separate $g$ into its regular and irregular parts as $g = g_\delta + \tilde g$ and denote its weighted version as
\begin{align*}
  G = \vint{v}^2 g, 
\qquad
  G = G_\delta + \tilde G,
\qquad
  G_\delta 
  = G \ast \eta_\delta
  = (\vint{v}^2 g) \ast \eta_\delta,
\end{align*}
where, as before, $\eta_\delta$ is the rescaled mollifier.
Decompose $\mathrm{Int}(I_{3})$ as
\begin{align} \label{decomp:I-3}
  \mathrm{Int}(I_{3})
&= \int_0^T \int_{\R^3} 
  \CalM^2 w \, \vint{v}^2 \nabla_v \cdot
  (A[w_0] \nabla_v g) \dv\dt \nn
\\
& \hspace{-1cm} = \int_0^T \int_{\R^3} 
  \CalM^2 w \, \nabla_v \cdot
  (A[w_0] \nabla_v G) \dv\dt
  - 2 \int_0^T \int_{\R^3} 
  (\CalM^2 w) \, v \cdot A[w_0] \nabla_v g  \dv\dt \nn
\\
& \qquad \,  
  - 2 \int_0^T \int_{\R^3} 
  (\CalM^2 w) \, \nabla_v \cdot (A[w_0] \, v \, g)  \dv\dt  \nn
\\
& \hspace{-1cm} = -\int_0^T \int_{\R^3} 
  \nabla_v \CalM(\CalM w) \,  \cdot
  (A[w_0] \nabla_v G_\delta) \dv\dt
  + \int_0^T \int_{\R^3}
    \nabla^2_v \CalM^2 w \,  : A[w_0] \, \tilde G \dv \dt \nn
\\
& + \int_0^T \int_{\R^3}
    \nabla_v \CalM^2 w \cdot (\nabla_v \cdot A[w_0]) \, \tilde G \dv \dt
  + 2 \int_0^T \int_{\R^3} 
  (\CalM^2 w) \, \nabla_v \cdot (v \cdot A[w_0]) \, g  \dv\dt \nn
\\
& \qquad \,  
  - 4 \int_0^T \int_{\R^3} 
  (\CalM^2 w) \, \nabla_v \cdot (A[w_0] \, v \, g)  \dv\dt   \nn
\\
& \hspace{-1cm}=: \mathrm{Int}(I_{3, 1})
   + \mathrm{Int}(I_{3, 2})
   + \mathrm{Int}(I_{3, 3})
   + \mathrm{Int}(I_{3, 4})
   + \mathrm{Int}(I_{3, 5}).
\end{align}
By~\eqref{reform:A-w} in Propositions~\ref{prop:A-divA-w},
\begin{align*}
  \abs{\mathrm{Int}(I_{3, 1})}
& \leq 
 \abs{\int_0^T \int_{\R^3} 
 \vpran{\nabla_v \CalM(\CalM w)}
 \vpran{A[\CalM w_0]} \nabla_v G_\delta \dv \dt }
\\
& \quad \,
 + \abs{\int_0^T \int_{\R^3} 
 \vpran{\nabla_v \CalM(\CalM w)}
 \vpran{A[\Delta \CalM w_0]} \nabla_v G_\delta \dv \dt}
=: 
     \mathrm{Int}(I_{3, 1}^{(1)})
     + \mathrm{Int}(I_{3, 1}^{(2)}).
\end{align*}
Applying~\eqref{bound:A-Mw0} in Proposition~\ref{prop:A-divA-w}, we have
\begin{align} \label{bound:I-311}
  \abs{\mathrm{Int}(I_{3, 1}^{(1)})}
&\leq
  \int_0^T \norm{\CalM w}_{L^2(\R^3)} 
  \vpran{\int_{\R^3} 
 \abs{\nabla_v \CalM(\CalM w)}
  |\vint{v} \nabla_v G_\delta| \dv} \dt  \nn
\\
& \leq
  \int_0^T \norm{\CalM w}_{L^2(\R^3)}
  \norm{\nabla_v \CalM (\CalM w)}_{L^2(\R^3)}
  \norm{\vint{v} \nabla_v G_\delta}_{L^2(\R^3)} \dt  \nn
\\
& \leq  
 \vpran{\Sup_{[0, T]} \norm{\vint{v} \nabla_v G_\delta}_{L^2(\R^3)}}
 \int_0^T
 \norm{\CalM w}_{L^2(\R^3)}^2 \dt  \nn
\\
& \leq
  \frac{C}{\delta^{3/2}} \vpran{\Sup_{[0, T]}\norm{\vint{\cdot}^3 g}_{L^{\frac{3}{2}}}}
  \int_0^T
 \norm{\CalM w}_{L^2(\R^3)}^2 \dt,
\end{align}
where we have applied Young's inequality to the $G_\delta$-term, which is bounded by
\begin{align*}
  \norm{\vint{v} \nabla_v G_\delta}_{L^2(\R^3)}
&\leq C \norm{(\vint{\cdot}^3  g) \ast \nabla \eta_\delta}_{L^2(\R^3)}
= \frac{1}{\delta}
  \norm{(\vint{\cdot}^3 g) \ast   \vpran{\nabla\eta_0}_\delta}_{L^2(\R^3)}
\\
&\leq
  \frac{1}{\delta}
  \norm{\vint{\cdot}^3 g}_{L^{\frac{3}{2}}}
  \norm{(\nabla \eta_0)_\delta}_{L^{\frac{6}{5}}}
\leq 
  \frac{1}{\sqrt{\delta}}
  \norm{\vint{\cdot}^3 g}_{L^{\frac{3}{2}}}.
\end{align*}
To bound the second term $\mathrm{Int}(I_{3, 1}^{(2)})$, we apply~\eqref{bound:A-Delta-Mw0} in Proposition~\ref{prop:A-divA-w}. This gives
\begin{align} \label{bound:I-313-1}
  \abs{ \mathrm{Int}(I_{3, 1}^{(2)})}
&\leq \int_0^T \int_{\R^3} 
 \abs{\nabla_v \CalM(\CalM w)}
 \abs{A[\Delta \CalM w]} |\nabla_v G_\delta| \dv \dt \nn
\\
& \leq
 \int_0^T
 \norm{\nabla_v \CalM(\CalM w)}_{L^2}
 \norm{A[\Delta \CalM w]}_{L^6}
 \norm{\nabla_v G_\delta}_{L^3} \dt \nn
\\
& \leq 
 C \vpran{\Sup_{[0, T]} \norm{\nabla_v G_\delta}_{L^3(\R^3)}}
 \int_0^T
  \norm{\CalM w}_{L^2}
  \vpran{\norm{\vint{v}^{-\frac{3}{2}} \nabla_v \CalM w}_{L^2}
  + \norm{\CalM w}_{L^2}} \dt \nn
\\
& \leq \epsilon_0 \int_0^T
 \norm{\vint{v}^{-\frac{3}{2}} \nabla_v\CalM w}^2_{L^2}
 \dt
 + \frac{C_{\epsilon_0}}{\delta^{5}}
 \int_0^T \norm{\CalM w}_{L^2}^2 \dt,
\end{align}
in which we have applied Young's inequality to get
\begin{align*}
  \norm{\nabla_v G_\delta}_{L^3(\R^3)}
\leq
 \frac{1}{\delta}
 \norm{G * (\nabla \eta_0)_\delta}_{L^3}
\leq
 \frac{C}{\delta}
 \norm{G}_{L^{3/2}}
 \norm{(\nabla \eta_0)_\delta}_{L^{3/2}}
\leq
  \frac{C}{\delta^{5/2}} \norm{\vint{\cdot}^2 g}_{L^{3/2}}.
\end{align*}
Together with \eqref{bound:I-311}, we obtain that 
\begin{align} \label{bound:I-31}
 \abs{ \mathrm{Int}(I_{3, 1}) }
\leq
  \epsilon_0 \int_0^T
 \norm{\vint{v}^{-\frac{3}{2}} \nabla_v\CalM w}^2_{L^2}
 \dt
 + \frac{C_{\epsilon_0}}{\delta^{5}}
 \int_0^T \norm{\CalM w}_{L^2}^2 \dt.
\end{align}

Similarly, we bound $\mathrm{Int}(I_{3, 2})$ as follows:
\begin{align*}
 \abs{ \mathrm{Int}(I_{3, 2})}
&\leq
  \int_0^T \!\! \int_{\R^3}
    \abs{\nabla^2_v \CalM^2 w} \abs{A[\CalM w_0]} \abs{\tilde G} \dv \dt
  + \int_0^T \!\! \int_{\R^3}
    \abs{\nabla^2_v \CalM^2 w} \abs{A[\Delta \CalM w_0]} \abs{\tilde G} \dv \dt.
\end{align*}
By Sobolev embedding $W^{\frac{1}{2}, 2} (\R^3) \hookrightarrow L^3(\R^3)$ we get
\begin{align*}
  \int_0^T \int_{\R^3}
    \abs{\nabla^2_v \CalM^2 w} \abs{A[\CalM w_0]} \abs{\tilde G} \dv \dt
&\leq
  \int_0^T 
  \norm{\CalM w}_{L^2(\R^3)}
  \int_{\R^3}
    \abs{\nabla^2_v \CalM^2 w} \abs{\vint{v}\tilde G} \dv \dt
\\
& \hspace{-2.5cm} \leq
  \int_0^T 
  \norm{\CalM w}_{L^2(\R^3)}
    \norm{\vint{v}^{-\frac{3}{2}} \nabla^2_v \CalM^2 w}_{L^3} \norm{\vint{v}^{\frac{5}{2}}\tilde G}_{L^{\frac{3}{2}}} \dt
\\
& \hspace{-2.5cm} \leq
  C \int_0^T 
  \norm{\CalM w}_{L^2(\R^3)}
    \norm{\nabla_v \CalM \vpran{\vint{v}^{-\frac{3}{2}} \nabla_v \CalM w}}_{L^3} \norm{\vint{v}^{\frac{5}{2}}\tilde G}_{L^{\frac{3}{2}}} \dt
\\
& \hspace{-2.5cm} \leq
  C \int_0^T 
  \norm{\CalM w}_{L^2(\R^3)}
    \norm{\nabla_v \CalM^{3/4} \vpran{\vint{v}^{-\frac{3}{2}} \nabla_v \CalM w}}_{L^3} \norm{\vint{v}^{\frac{5}{2}}\tilde G}_{L^{\frac{3}{2}}} \dt
\\
& \hspace{-2.5cm} \leq
  C \int_0^T 
  \norm{\CalM w}_{L^2(\R^3)}
    \norm{\vint{v}^{-\frac{3}{2}} \nabla_v \CalM w}_{L^2} \norm{\vint{v}^{\frac{5}{2}}\tilde G}_{L^{\frac{3}{2}}} \dt
\\
& \hspace{-2.5cm} \leq
 C \vpran{\Sup_{[0, T]} \norm{\vint{v}^{\frac{9}{2}}\tilde g}_{L^{\frac{3}{2}}}}
 \int_0^T 
  \norm{\CalM w}_{L^2(\R^3)}
    \norm{\vint{v}^{-\frac{3}{2}} \nabla_v \CalM w}_{L^2}  \dt
\\
& \hspace{-2.5cm} \leq  
  \epsilon_0 \int_0^T 
  \norm{\vint{v}^{-\frac{3}{2}} \nabla_v \CalM w}^2_{L^2} \dt
  + C_{\epsilon_0} \int_0^T 
  \norm{\CalM w}^2_{L^2(\R^3)} \dt.
\end{align*}
Furthermore, by Proposition~\ref{prop:A-divA-w}, 
\begin{align*}
& \quad \,
  \int_0^T \int_{\R^3}
    \abs{\nabla^2_v \CalM^2 w} \abs{A[\Delta \CalM w_0]} \abs{\tilde G} \dv \dt
\\
& \leq
\int_0^T
 \norm{\vint{v}^{-\frac{3}{2}}\nabla \CalM (\nabla \CalM w)}_{L^6}
 \norm{A[\Delta \CalM w_0]}_{L^6}
 \norm{\vint{v}^{3/2} \tilde G}_{L^{\frac{3}{2}}} \dt
\\
& \leq 
 C \vpran{\Sup_{[0, T]} \norm{\vint{v}^{3/2} \tilde G}_{L^{\frac{3}{2}}}}
 \int_0^T 
 \vpran{\norm{\vint{v}^{-\frac{3}{2}} \nabla_v \CalM w}^2_{L^2} + \norm{\CalM w}^2_{L^2}} \dt
\\
& \leq
 \epsilon_0 
 \int_0^T 
 \norm{\vint{v}^{-\frac{3}{2}} \nabla_v \CalM w}^2_{L^2} \dt
 + C \int_0^T \norm{\CalM w}^2_{L^2} \dt,
\end{align*}
by taking $\delta$ small enough. 
Overall, we have
\begin{align} \label{bound:I-32}
  \abs{\mathrm{Int}(I_{3, 2})}
&\leq
  \epsilon_0 \int_0^T 
  \norm{\vint{v}^{-\frac{3}{2}} \nabla_v \CalM w}^2_{L^2} \dt
  + C_{\epsilon_0} \int_0^T 
  \norm{\CalM w}^2_{L^2(\R^3)} \dt.
\end{align}
By~\eqref{reform:div-A-w} in Proposition~\ref{prop:A-divA-w}, we bound $\mathrm{Int}(I_{3, 3})$ as follows. 
\begin{align*}
  \abs{\mathrm{Int}(I_{3, 3})}
&= \abs{\int_0^T \int_{\R^3}
    \nabla_v \CalM^2 w \cdot (\nabla_v \cdot A[w]) \, \tilde G \dv \dt}
\\
& \leq 
  \frac{1}{4 \pi} \int_0^T \int_{\R^3}
    \abs{\nabla_v \CalM^2 w} 
    \abs{A_1} \abs{\tilde G} \dv \dt
    + \int_0^T \int_{\R^3}
    \abs{\nabla_v \CalM^2 w}
    \abs{\nabla_v \CalM w_0}
    \abs{\tilde G} \dv \dt
\\
& 
=: \mathrm{Int}(I_{3, 3}^{(1)})
   + \mathrm{Int}(I_{3, 3}^{(2)}), 
\end{align*}
where by the bound of $A_1$ in~\eqref{reform:div-A-w},  
\begin{align*}
 \abs{ \mathrm{Int}(I_{3, 3}^{(1)})}
& \leq 
  C \vpran{\Sup_{[0, T]} \norm{\tilde G}_{L^{\frac{3}{2}}}}
  \int_0^T 
  \norm{\nabla_v \CalM^2 w}_{L^6}
  \norm{\CalM w_0}_{L^2} \dt
 \leq
  C \int_0^T 
  \norm{\CalM w}^2_{L^2(\R^3)} \dt.
\end{align*}
Finally, by H\"{o}lder's inequality and the relation $w = \vint{v}^2 w_0$, 
\begin{align*}
  \abs{\mathrm{Int}(I_{3, 3}^{(2)})}
&\leq
  \int_0^T
  \norm{\nabla_v \CalM^2 w}_{L^6}
  \norm{\vint{v}^{-3/2}\nabla_v \CalM w}_{L^2}
  \norm{\tilde G}_{L^3} \dt
\\
& \quad \,
  + \int_0^T
  \norm{\nabla_v \CalM^2 w}_{L^6}
  \norm{\CalM w}_{L^2}
  \norm{\tilde G}_{L^3} \dt
\end{align*}
Note that by H\"{o}lder's inequality again,
\begin{align*}
  \norm{\tilde G}_{L^3}
&\leq
  \norm{\tilde G}_{L^{9/2}}^{3/4}
  \norm{\tilde G}_{L^{3/2}}^{1/4}
\\
&\leq
 C \norm{\tilde G}_{L^{3/2}}^{1/4}
 \norm{\vint{v}^{-\frac{3}{2}}\nabla (\vint{v}^4 g)^{3/4}}_{L^2}
 + C \norm{\tilde G}_{L^{3/2}},
\end{align*}
where $C$ is independent of $\delta$.
Therefore,
\begin{align*}
  \mathrm{Int}(I_{3, 3}^{(2)})
\leq
 &\, C \vpran{\Sup_{[0, T]} \norm{\tilde G}_{L^{\frac{3}{2}}}}^{\frac{1}{4}}
  \vpran{\Sup_{[0, T]}
  \norm{\CalM w}_{L^2}}
\\
& \hspace{1cm} \times
  \int_0^T   \vpran{\norm{\vint{v}^{-3/2}\nabla_v \CalM w}_{L^2} + \norm{\CalM w}_{L^2}}
   \norm{\vint{v}^{-\frac{3}{2}}\nabla (\vint{v}^4 g)^{3/4}}_{L^2} \dt
\\
& \quad \,
  + C \vpran{\Sup_{[0, T]} \norm{\tilde G}_{L^{\frac{3}{2}}}}
  \int_0^T
  \norm{\CalM w}_{L^2}
  \norm{\vint{v}^{-\frac{3}{2}}\nabla \CalM w}_{L^2} \dt
\\
& \hspace{-2cm}\leq
  C \vpran{\Sup_{[0, T]} \norm{\tilde G}_{L^{\frac{3}{2}}}}^{\frac{1}{2}}
  \vpran{\int_0^T \norm{\vint{v}^{-\frac{3}{2}}\nabla (\vint{v}^5 g)^{3/4}}^2_{L^2} \dt}^{1/2}
  \vpran{\Sup_{[0, T]}
  \norm{\CalM w}^2_{L^2}}
\\
& \hspace{-1.5cm} 
  + C \vpran{\Sup_{[0, T]} \norm{\tilde G}_{L^{\frac{3}{2}}}}^{\frac{1}{2}}
  \vpran{\int_0^T \norm{\vint{v}^{-\frac{3}{2}}\nabla (\vint{v}^5 g)^{3/4}}^2_{L^2} \dt}^{1/2}
  \int_0^T
  \norm{\vint{v}^{-3/2}\nabla_v \CalM w}^2_{L^2} \dt
\\
& \hspace{-1.5cm} 
  + C \vpran{\Sup_{[0, T]} \norm{\tilde G}_{L^{\frac{3}{2}}}}
  \int_0^T \norm{\CalM w}^2_{L^2} \dt.
\end{align*}
By taking $\delta$ small enough, we have
\begin{align*}
   \abs{ \mathrm{Int}(I_{3, 3}^{(2)})}
\leq
  \epsilon_0 \Sup_{[0, T]}
  \norm{\CalM w}^2_{L^2}
  + \epsilon_0 
  \int_0^T
  \norm{\vint{v}^{-3/2}\nabla_v \CalM w}^2_{L^2} \dt
  + C \int_0^T
    \norm{\CalM w}^2_{L^2} \dt.
\end{align*}
Therefore, we obtain the bound for $\mathrm{Int}(I_{3, 3})$ as
\begin{align} \label{bound:I-33}
  \abs{\mathrm{Int}(I_{3, 3})}
&\leq
  \epsilon_0 \, \Sup_{[0, T]}
  \norm{\CalM w}^2_{L^2}
  + \epsilon_0  \!
  \int_0^T
  \norm{\vint{v}^{-3/2} \nabla_v \CalM w}^2_{L^2} \dt  
  + C_{\epsilon_0} \int_0^T 
  \norm{\CalM w}^2_{L^2} \dt.
\end{align}

The last two terms, $\mathrm{Int}(I_{3,4})$ and $\mathrm{Int}(I_{3,5})$, are lower order terms that can be bounded directly. More specifically, 
\begin{align} \label{bound:I-3-4}
 \abs{\mathrm{Int}(I_{3,4})}
&\leq
 2 \int_0^T \int_{\R^3} 
  \abs{\CalM^2 w} \, \abs{\nabla_v \cdot (v \cdot A[w_0]} \, |g|  \dv\dt  \nn
\\
& \hspace{-1cm} \leq
 C \int_0^T \int_{\R^3} 
  \abs{\CalM^2 w} \, \abs{A[w_0]} \, |g|  \dv\dt
 + C \int_0^T \int_{\R^3} 
  \abs{\CalM^2 w} \, \abs{\nabla_v \cdot A[w_0]} \, |v g|  \dv\dt \nn
\\
& \hspace{-1cm} =: \mathrm{Int}(I_{3,4}^{(1)})
 + \mathrm{Int}(I_{3,4}^{(2)}).
\end{align}
Compared with $\mathrm{Int}(I_{3,2})$, $\mathrm{Int}(I_{3,4}^{(1)})$ has a similar and more regular structure, with 
$\nabla_v^2 \CalM^2 w$ replaced by $\CalM^2 w$. Therefore, it is expected that $\mathrm{Int}(I_{3,4}^{(1)})$ satisfies a similar bound as $\mathrm{Int}(I_{3,2})$ without the need to separate $g$ into its smooth and non-smooth parts. Indeed, by Proposition~\ref{prop:A-divA-w},
\begin{align} \label{decomp:I-34-1}
 \abs{ \mathrm{Int}(I_{3,4}^{(1)})}
&\leq
  C \! \int_0^T \!\! \int_{\R^3} 
  \abs{\CalM^2 w} \, \abs{A[\CalM w_0]} \, |g|  \dv\dt  
  + C \! \int_0^T \!\! \int_{\R^3} 
  \abs{\CalM^2 w} \, \abs{A[\Delta \CalM w_0]} \, |g|  \dv\dt,
\end{align}
where, by~\eqref{reform:A-w} in Proposition~\ref{prop:A-divA-w} and Sobolev embedding $W^{\frac{1}{2}, 2} (\R^3) \hookrightarrow L^3(\R^3)$,
\begin{align*} 
  \int_0^T \int_{\R^3} 
  \abs{\CalM^2 w} \, \abs{A[\CalM w_0]} \, |g|  \dv\dt
&\leq
  C \int_0^T \norm{\CalM w}_{L^2} \vpran{\int_{\R^3} 
  \abs{\CalM^2 w}  |v g|  \dv}\dt \nn
\\
& \hspace{-1cm} \leq
  C \vpran{\Sup_{[0, T]} \norm{v g}_{L^{3/2}}}\int_0^T \norm{\CalM w}_{L^2}  
  \norm{\CalM^2 w}_{L^3} \dt \nn
\\
& \hspace{-1cm} \leq
  C \vpran{\Sup_{[0, T]} \norm{v g}_{L^{3/2}}}\int_0^T \norm{\CalM w}_{L^2}  
  \norm{\CalM^{\frac{1}{4}} \CalM w}_{L^3} \dt \nn
\\
& \hspace{-1cm} \leq 
  C \vpran{\Sup_{[0, T]} \norm{v g}_{L^{3/2}}}\int_0^T \norm{\CalM w}_{L^2}^2 \dt,
\end{align*}
and by~\eqref{reform:div-A-w} in Proposition~\ref{prop:A-divA-w} together with Sobolev embedding $W^{1, 2} (\R^3) \hookrightarrow L^6(\R^3)$,
\begin{align*}
& \quad \,
  \int_0^T \int_{\R^3} 
  \abs{\CalM^2 w} \, \abs{A[\Delta \CalM w_0]} \, |g|  \dv\dt
\\
& \leq 
 \vpran{\Sup_{[0, T]} \norm{g}_{L^{3/2}}}
 \int_0^T
 \norm{\CalM^2 w}_{L^6}
 \norm{A[\Delta \CalM w_0]}_{L^6} \dt
\\
& \leq 
 \vpran{\Sup_{[0, T]} \norm{g}_{L^{3/2}}}
 \int_0^T
 \norm{\CalM^{1/2} \CalM w}_{L^6}
 \norm{\vint{\cdot}^{-3/2}\nabla \CalM w}_{L^2} \dt
\\
& \leq 
 \vpran{\Sup_{[0, T]} \norm{g}_{L^{3/2}}}
 \int_0^T
 \norm{\CalM w}_{L^2}
 \norm{\vint{\cdot}^{-3/2}\nabla \CalM w}_{L^2} \dt
\\
& \leq
 \epsilon_0 \int_0^T
 \norm{\vint{\cdot}^{-3/2}\nabla \CalM w}_{L^2}^2 \dt
 + C_{\epsilon_0} \int_0^T
   \norm{\CalM w}_{L^2}^2 \dt. 
\end{align*}
Therefore,
\begin{align} \label{bound:I-34-1}
  \abs{\mathrm{Int}(I_{3,4}^{(1)})}
&\leq
  \epsilon_0 \int_0^T
 \norm{\vint{\cdot}^{-3/2}\nabla \CalM w}_{L^2}^2 \dt
 + C_{\epsilon_0} \int_0^T
   \norm{\CalM w}_{L^2}^2 \dt.
\end{align}
The estimate for $\mathrm{Int}(I_{3,4}^{(2)})$ is similar to $\mathrm{Int}(I_{3,3})$. By~\eqref{reform:div-A-w} in Proposition~\ref{prop:A-divA-w},
\begin{align*}
\abs{  \mathrm{Int}(I_{3,4}^{(2)})}
& \leq 
  C \int_0^T \int_{\R^3} 
  \abs{\CalM^2 w} \, \abs{A_1} \, |v g|  \dv\dt
  + C \int_0^T \int_{\R^3} 
  \abs{\CalM^2 w} \, \abs{\nabla \CalM w_0} \, |v g|  \dv\dt,
\end{align*}
where, by~\eqref{bound:A-1}, the first term satisfies
\begin{align} \label{bound:lot-A1}
  \int_0^T \int_{\R^3} 
  \abs{\CalM^2 w} \, \abs{A_1} \, |v g|  \dv\dt
&\leq
 \vpran{\Sup_{[0, T]} \norm{v g}_{L^{3/2}}}
 \int_0^T
 \norm{\CalM^2 w}_{L^6} \norm{A_1}_{L^6} \dt \nn
\\
& \leq 
 C \int_0^T
   \norm{\CalM w}_{L^2}^2 \dt,
\end{align}
and 
\begin{align} \label{bound:lot-1}
 \int_0^T \int_{\R^3} 
  \abs{\CalM^2 w} \, \abs{\nabla \CalM w_0} \, |v g|  \dv\dt
&\leq
  \int_0^T \int_{\R^3} 
  \abs{\CalM^2 w} \, \abs{\nabla \CalM w_0} \, |(v g)_\delta|  \dv\dt  \nn
\\
& \quad \,
  + \int_0^T \int_{\R^3} 
  \abs{\CalM^2 w} \, \abs{\nabla \CalM w_0} \, |\widetilde {vg}|  \dv\dt,
\end{align}
where $\widetilde {vg} = v g - (v g)_{\delta} = v g - (v g) \ast \eta_\delta$. The first term on the right-hand side of~\eqref{bound:lot-1} satisfies
\begin{align} \label{bound:lot-2}
& \quad \,
  \int_0^T \int_{\R^3} 
  \abs{\CalM^2 w} \, \abs{\nabla \CalM w_0} \, |(v g)_\delta|  \dv\dt \nn
\\
&\leq
 \vpran{\Sup_{[0, T]} \norm{(v g)_\delta}_{L^\infty}}
 \int_0^T 
 \norm{\CalM^2 w}_{L^2}
 \norm{\nabla \CalM w_0}_{L^2} \dt \nn
\\
& \leq 
 C_\delta \vpran{\Sup_{[0, T]} \norm{v g}_{L^{3/2}}}
 \int_0^T 
 \norm{\CalM w}_{L^2}
 \norm{\vint{v}^{-3/2} \nabla \CalM w}_{L^2} \dt \nn
\\
& \quad \,
 + C_\delta \vpran{\Sup_{[0, T]} \norm{v g}_{L^{3/2}}}
 \int_0^T 
 \norm{\CalM w}_{L^2}^2 \dt \nn
\\
& \leq
 \epsilon_0 \int_0^T
 \norm{\vint{v}^{-3/2} \nabla \CalM w}_{L^2}^2 \dt
 + C_{\epsilon_0} 
 \int_0^T \norm{\CalM w}_{L^2}^2 \dt.
\end{align}
The second term on the right-hand side of~\eqref{bound:lot-1} satisfies
\begin{align*}
& \quad \,
  \int_0^T \int_{\R^3} 
  \abs{\CalM^2 w} \, \abs{\nabla \CalM w_0} \, |\widetilde {vg}|  \dv\dt
\\
& \leq
  \int_0^T \int_{\R^3} 
  \norm{\CalM^2 w}_{L^6} \, \norm{\nabla \CalM w_0}_{L^2} \, \norm{\widetilde {vg}}_{L^3}  \dv\dt,
\end{align*}
which has a similar structure as $\mathrm{Int}(I_{3, 3}^{(2)})$. The only difference is that we replace $\norm{\nabla \CalM^2 w}_{L^6}$ in $\mathrm{Int}(I_{3, 3}^{(2)})$ by $\norm{\CalM^2 w}_{L^6}$, and both of them are bounded by $C \norm{\CalM w}_{L^2}$. Repeating the estimates for $\mathrm{Int}(I_{3, 3}^{(2)})$, we obtain
\begin{align} \label{bound:lot-3}
  \int_0^T \int_{\R^3} 
  \abs{\CalM^2 w} \, \abs{\nabla \CalM w_0} \, |\widetilde {vg}|  \dv\dt
&\leq
  \epsilon_0 \Sup_{[0, T]}
  \norm{\CalM w}^2_{L^2}
  + \epsilon_0 
  \int_0^T
  \norm{\vint{v}^{-3/2}\nabla_v \CalM w}^2_{L^2} \dt \nn
\\
& \quad \,
  + C \int_0^T
    \norm{\CalM w}^2_{L^2} \dt. 
\end{align}
Combining~\eqref{bound:lot-3} with~\eqref{bound:lot-2} and~\eqref{bound:lot-A1}, we obtain that
\begin{align} \label{bound:I-34-2}
 \abs{ \mathrm{Int}(I_{3,4}^{(2)})}
&\leq
  \epsilon_0 \Sup_{[0, T]}
  \norm{\CalM w}^2_{L^2}
  + \epsilon_0 
  \int_0^T
  \norm{\vint{v}^{-3/2}\nabla_v \CalM w}^2_{L^2} \dt 
  + C \int_0^T
    \norm{\CalM w}^2_{L^2} \dt.
\end{align}
Combining~\eqref{bound:I-34-2} with~\eqref{bound:I-34-1} gives
\begin{align}\label{bound:I-34}
  \abs{\mathrm{Int}(I_{3,4})}
&\leq
  \epsilon_0 \Sup_{[0, T]}
  \norm{\CalM w}^2_{L^2}
  + \epsilon_0 
  \int_0^T
  \norm{\vint{v}^{-3/2}\nabla_v \CalM w}^2_{L^2} \dt
  + C \int_0^T
    \norm{\CalM w}^2_{L^2} \dt.
 \end{align}

To bound the last term $\mathrm{Int}(I_{3,5})$ in~\eqref{decomp:I-3}, we integrate by parts and apply the decomposition in Proposition~\ref{prop:A-divA-w}:
\begin{align} \label{decomp:I-35}
 \abs{\mathrm{Int}(I_{3,5})}
&= 4\abs{\int_0^T \int_{\R^3} 
  \nabla_v \CalM^2 w \cdot (A[w_0] \, v \, g)  \dv\dt} \nn
\\
&\leq
 C \int_0^T \int_{\R^3} 
  \abs{\nabla_v \CalM^2 w} \abs{A[w_0]} \abs{v g}  \dv\dt \nn
\\
& \leq
 C \int_0^T \int_{\R^3} 
  \abs{\nabla_v \CalM^2 w} \abs{A[\CalM w_0]} \abs{v g}  \dv\dt \nn
\\
& \quad \,
  + C \int_0^T \int_{\R^3} 
  \abs{\nabla_v \CalM^2 w} \abs{A[\Delta \CalM w_0]} \abs{v g}  \dv\dt,
\end{align}
where by~\eqref{bound:A-Mw0}, the first term on the right-hand side of~\eqref{decomp:I-35} satisfies
\begin{align*}
  \int_0^T \int_{\R^3} 
  \abs{\nabla_v \CalM^2 w} \abs{A[\CalM w_0]} \abs{v g}  \dv\dt
&\leq 
 \int_0^T
 \norm{\CalM w}_{L^2}
 \int_{\R^3}
 \abs{\nabla_v \CalM (\CalM w)}
 \abs{\vint{v}^2 g} \dv \dt
\\
& \hspace{-1cm} \leq
 \int_0^T
 \norm{\CalM w}_{L^2}
 \norm{\nabla_v \CalM^{3/4} (\CalM w)}_{L^3}
 \norm{\vint{v}^2 g}_{L^{3/2}} \dt
\\
& \hspace{-1cm} \leq
 C \vpran{\Sup_{[0, T]} \norm{\vint{v}^2 g}_{L^{3/2}}}
 \int_0^T \norm{\CalM w}_{L^2}^2 \dt,
\end{align*}
and the second term is bounded similarly to the second term for $\mathrm{Int}(I_{3,4}^{(1)})$ in~\eqref{decomp:I-34-1}, with $\norm{\CalM^2 w}_{L^6}$ in $\mathrm{Int}(I_{3,4}^{(1)})$ replaced by $\norm{\nabla \CalM^2 w}_{L^6}$ here. Therefore, we have
\begin{align} \label{bound:I-35}
    \abs{\mathrm{Int}(I_{3,5})}
\leq
  \epsilon_0 \int_0^T
 \norm{\vint{\cdot}^{-3/2}\nabla \CalM w}_{L^2}^2 \dt
 + C_{\epsilon_0} \int_0^T
   \norm{\CalM w}_{L^2}^2 \dt.
\end{align}

Combining~\eqref{bound:I-31}, \eqref{bound:I-32}, \eqref{bound:I-33}, \eqref{bound:I-34}, and \eqref{bound:I-35}, we have
\begin{align*}
  \abs{\mathrm{Int}(I_3)}
\leq 
    \epsilon_0 \, \Sup_{[0, T]}
  \norm{\CalM w}^2_{L^2}
  + \epsilon_0 
  \int_0^T
  \norm{\vint{v}^{-3/2} \nabla_v \CalM w}^2_{L^2} \dt  \nn
  + C_{\epsilon_0} \int_0^T 
  \norm{\CalM w}^2_{L^2(\R^3)} \dt.
\end{align*}
To bound $\mathrm{Int}(I_4)$, we first recall that $G = \vint{v}^2 g$ and 
\begin{align*}
  \mathrm{Int}(I_4)
&= - \int_0^T \int_{\R^3} 
  \CalM^2 w \, \vint{v}^2 \nabla_v \cdot
  (g \nabla_v a[w_0]) \dv\dt
\\
&= - \int_0^T \int_{\R^3} 
  \CalM^2 w \, \nabla_v \cdot
  (G \, \nabla_v a[w_0]) \dv\dt
  + 2 \int_0^T \int_{\R^3} 
  \CalM^2 w \, (v g)
  \cdot \nabla_v a[w_0] \dv\dt
\\
& = \int_0^T \int_{\R^3} 
  \nabla_v \CalM^2 w \,  \cdot
  (\nabla_v a[w_0]) \, G \dv\dt
  + 2 \int_0^T \int_{\R^3} 
  \CalM^2 w \, (v g)
  \cdot \nabla_v a[w_0] \dv\dt
\\
&= \int_0^T \int_{\R^3} 
  \nabla_v \CalM^2 w \,  \cdot
  (\nabla_v a[w]) \, G_\delta \dv\dt
  + \int_0^T \int_{\R^3} 
  \nabla_v \CalM^2 w \,  \cdot
  (\nabla_v \cdot A[w]) \tilde G \dv\dt
\\
& \quad \,
 + 2 \int_0^T \int_{\R^3} 
  \CalM^2 w \, (v g)
  \cdot \nabla_v a[w_0] \dv\dt
=: \mathrm{Int}(I_{4,1})
  + \mathrm{Int}(I_{4,2})
  + \mathrm{Int}(I_{4,3}),
\end{align*}
where $\mathrm{Int}(I_{4,2})$ is the same as $\mathrm{Int}(I_{3,3})$, and $\mathrm{Int}(I_{4,3})$ is upper-bounded in the same way as $\mathrm{Int}(I_{3,4}^{(2)})$. By Proposition~\ref{prop:A-divA-w}, $\mathrm{Int}(I_{4,1})$ satisfies
\begin{align*}
& \quad \,
  \abs{\int_0^T \int_{\R^3} 
  \nabla_v \CalM^2 w \,  \cdot
  (\nabla_v a[w_0]) \, G_\delta \dv\dt}
\\
& \leq 
  C \vpran{\Sup_{[0, T]}
    \norm{G_\delta}_{L^{3/2}}}
  \int_0^T
   \norm{\nabla_v \CalM^2 w}_{L^6}
   \norm{A_1}_{L^6} \dt
\\
& \quad \,
  + C \int_0^T
   \norm{\nabla_v \CalM^2 w}_{L^6}
   \norm{\nabla_v \CalM w_0}_{L^2}
   \norm{G_\delta}_{L^3} \dt
\\
& \leq
  C\vpran{\Sup_{[0, T]}
    \norm{\vint{\cdot}^2 g}_{L^{3/2}}}
  \int_0^T
   \norm{\CalM w}^2_{L^2} \dt
\\
& \quad \, 
  + C \vpran{\Sup_{[0, T]}
    \norm{G_\delta}_{L^3}}
    \vpran{\int_0^T \norm{\CalM w}^2_{L^2} \dt}^{1/2}
    \vpran{\int_0^T \norm{\vint{v}^{-3/2} \nabla_v \CalM w}^2_{L^2} \dt}^{1/2}
\\
& \leq
  \frac{C_{\epsilon_0}}{\delta^3}
  \vpran{\Sup_{[0, T]}\norm{\vint{\cdot}^2 g}^2_{L^{3/2}}}
  \int_0^T
   \norm{\CalM w}^2_{L^2} \dt
  + \epsilon_0
    \int_0^T \norm{\vint{v}^{-3/2} \nabla_v \CalM w}^2_{L^2} \dt
\\
& \leq
  \frac{C_{\epsilon_0}}{\delta^3}
  \vpran{\Sup_{[0, T]}\norm{\vint{\cdot}^2 g}^2_{L^{3/2}}}
  \int_0^T
   \norm{\CalM w}^2_{L^2} \dt
  + \epsilon_0
    \int_0^T \norm{\vint{v}^{-3/2} \nabla_v \CalM w}^2_{L^2} \dt.
\end{align*}
Thus, $\mathrm{Int}(I_4)$ satisfies a similar bound as $\mathrm{Int}(I_3)$:
\begin{align*}
  \abs{\mathrm{Int}(I_4)}
&\leq
 \epsilon_0
    \int_0^T \norm{\vint{v}^{-3/2} \nabla_v \CalM w}^2_{L^2} \dt
 + \epsilon_0 \, \Sup_{[0, T]} \norm{\CalM w}^2_{L^2}
 + C_{\epsilon_0}
   \int_0^T 
   \norm{\CalM w}^2_{L^2} \dt. \qedhere
\end{align*}
\end{proof}

\medskip

\subsection{Estimates of $\mathrm{Int}(I_5)$ and $\mathrm{Int}(I_6)$} The last two terms, $\mathrm{Int}(I_5)$ and $\mathrm{Int}(I_6)$, are of lower order, and they satisfy similar bounds as some of the terms in $I_1$ and $I_2$. We will sketch the proof of their bounds and omit some repetitive details. More precisely, we have

\begin{prop} \label{prop:I-6-I-7}
Let $\epsilon_0 > 0$ be an arbitrarily small number. Then there exists $C_{\epsilon_0} > 0$ such that 
\begin{align} \label{bound:I-5-6}
 \abs{ \mathrm{Int}(I_{5})}
  + \abs{\mathrm{Int}(I_{6})}
&\leq 
    \epsilon_0 \, \Sup_{[0, T]}
  \norm{\CalM w}^2_{L^2}
  + \epsilon_0 
  \int_0^T
  \norm{\vint{v}^{-3/2} \nabla_v \CalM w}^2_{L^2} \dt  \nn
\\
& \hspace{6cm}
  + C_{\epsilon_0} \int_0^T 
  \norm{\CalM w}^2_{L^2(\R^3)} \dt.
\end{align}  
\end{prop}

\begin{proof}
We start with $I_6$, which satisfies
\begin{align*}
  \mathrm{Int}(I_6)
= 2 \int_0^T \int_{\R^3} 
  \CalM^2 w \, (v \, w_0) \cdot \nabla_v a[f] \dv\dt.
\end{align*}
Note that $\mathrm{Int}(I_6)$ has a similar and more regular structure to $\mathrm{Int}(I_2)$, where $\nabla \CalM^2 w$ and $w$ in $\mathrm{Int}(I_2)$ are replaced by $\CalM^2 w$ and $v w_0 = \vint{v}^{-1} w$ respectively. Hence, $\mathrm{Int}(I_6)$ can be bounded in a similar way as $\mathrm{Int}(I_2)$, which gives the part of bound for $\mathrm{Int}(I_6)$ in~\eqref{bound:I-5-6}.

To bound $I_5$, we first rewrite
\begin{align*}
 R_1
&= -2v \cdot A [f] \nabla_v w_0
    - 2 \nabla_v \cdot (A[f] \, v \, w_0)
\\
&= - 4 \nabla_v \cdot (A[f] \, v \, w_0)
  + 2 \nabla_v \cdot (v \cdot A[f]) w_0
\\
&= - 4 \nabla_v \cdot (A[f] \, v \, w_0)
  + 2 \, a[f] \, w_0
  + 2 (v \cdot \nabla_v a[f]) w_0,
\end{align*}
where the last term is $R_2$ and we have used $Tr A[f] = a[f]$. Hence, 
\begin{align} \label{decomp:I-7}
  \mathrm{Int}(I_5)
&= 4 \int_0^T \int_{\R^3} 
  \nabla_v \CalM^2 w \, \cdot
  A[f] \, v \, w_0 \dv\dt \nn
  + 2 \int_0^T \int_{\R^3}
  \CalM^2 w \, a[f] \, w_0 \dv\dt
  + \mathrm{Int}(I_6)
\\
&=: \mathrm{Int}(I_{5,1})
    + \mathrm{Int}(I_{5,2})
    + \mathrm{Int}(I_6).
\end{align}
We only need to bound the first two terms. By the definition of $A$, we have
\begin{align*}
  A[f] \, v
= \frac{1}{8\pi} \int_{\mathbb{R}^3} \frac{P(v-v_1) }{|v-v_1|} v f(t,v_1)\dv_1
= \frac{1}{8\pi} \int_{\mathbb{R}^3} \frac{P(v-v_1) }{|v-v_1|} v_1 f(t,v_1)\dv_1
= A[v f],
\end{align*}
since $P(v - v_1) (v-v_1) = 0$. Therefore,
\begin{align*}
  \mathrm{Int}(I_{5,1})
&= 4 \int_0^T \int_{\R^3} 
  \nabla_v \CalM^2 w \, \cdot
  A[v f] \, w_0 \dv\dt
\\
&= 4 \int_0^T \int_{\R^3} 
  \nabla_v \CalM^2 w \, \cdot
  A[v f] \, (I - \Delta) \CalM w_0 \dv\dt
\\
&= 4 \int_0^T \int_{\R^3} 
  (I - \Delta) \vpran{\nabla_v \CalM^2 w \, \cdot
  A[v f]} \CalM w_0 \dv\dt
\\
&= 4 \int_0^T \int_{\R^3}
   \nabla_v \CalM w \cdot
  A[v f] \CalM w_0 \dv\dt
  -8  \int_0^T \int_{\R^3}
   \nabla_v^2 \CalM^2 w : \nabla_v A[v f] \, \CalM w_0 \dv\dt
\\
& \quad \,
  - 4\int_0^T \int_{\R^3}
   \nabla_v \CalM^2 w \cdot
   (\Delta_v A[v f]) \CalM w_0 \dv\dt
=: \mathrm{Int}(I_{5,1}^{(1)})
   + \mathrm{Int}(I_{5,1}^{(2)})
   + \mathrm{Int}(I_{5,1}^{(3)}).
\end{align*}
The term $\mathrm{Int}(I_{5,1}^{(2)})$ has a similar structure to $I_{1,2}^{(1)}$, with $\nabla_v^3 \CalM^2 w$ and $\CalM w$ in $I_{1,2}^{(1)}$ replaced by the more regular $\nabla_v \CalM^2 w$ and $\CalM w_0$, respectively, and $f$ replaced by $vf$. Moreover, $\mathrm{Int}(I_{5,1}^{(3)})$ has a similar structure to the second term in~\eqref{decomp:I-12-interm}, with $\CalM w$ in~\eqref{decomp:I-12-interm} replaced by $\CalM w_0$, $f$ replaced by $vf$ and $\nabla_v^2 \CalM^2 w$ replaced by the more regular $\nabla_v \CalM^2 w$. Therefore, they have similar bounds to $I_{1,2}$. We focus on $\mathrm{Int}(I_{5,1}^{(1)})$, which can be rewritten as
\begin{align*}
  \mathrm{Int}(I_{5,1}^{(1)})
&= 4 \int_0^T \int_{\R^3}
   \nabla_v \CalM w \cdot
  A[(v f)_\delta] \CalM w_0 \dv\dt
  + 4 \int_0^T \int_{\R^3}
   \nabla_v \CalM w \cdot
  A[\widetilde{v f}] \CalM w_0 \dv\dt
\\
&=: \mathrm{Int}(I_{5,1}^{(1,1)})
    + \mathrm{Int}(I_{5,1}^{(1,2)}),
\end{align*}
where we have decomposed $vf$ into its smooth and less regular parts: $(v f)_\delta = (v f) \ast \eta_\delta$ and $\widetilde{v f} = v f - (v f)_\delta$. Using the smoothness, we can bound $\mathrm{Int}(I_{5,1}^{(1,1)})$ directly by
\begin{align*}
  \abs{\mathrm{Int}(I_{5,1}^{(1,1)})}
&\leq
  C \vpran{\Sup_{[0, T]} \norm{A[(v f)_\delta]}_{L^\infty}}
  \vpran{\int_0^T \norm{\vint{v}^{-3/2} \nabla_v \CalM w}^2_{L^2} \dt}^{1/2}
\\
& \hspace{7cm} \times
  \vpran{\int_0^T \norm{\vint{v}^{3/2}\CalM w_0}^2_{L^2} \dt}^{1/2}.
\end{align*}
Note that for any $h \in L^1 \cap L^\infty$, we have
\begin{align*}
  |A[h]|
&\leq
  C \int_{\R^3} \frac{|h(v_1)| }{|v-v_1|} \dv_1
\\
&= C \int_{|v-v_1|\leq 1} \frac{|h(v_1)| }{|v-v_1|}  \dv_1
+ C \int_{|v-v_1|\geq 1} \frac{|h(v_1)| }{|v-v_1|} \dv_1
\leq
  C \norm{h}_{L^1\cap L^\infty}.
\end{align*}
Therefore,
\begin{align*}
  \abs{\mathrm{Int}(I_{5,1}^{(1,1)})}
&\leq
  C \vpran{\Sup_{[0, T]}(\norm{(vf)_\delta}_{L^1 \cap L^\infty})}
  \vpran{\int_0^T \norm{\vint{v}^{-3/2} \nabla_v \CalM w}^2_{L^2} \dt}^{1/2}
\\
& \hspace{7cm} \times
  \vpran{\int_0^T \norm{\CalM w}_{L^2}^2 \dt}^{1/2}
\\
& \hspace{-1cm}\leq 
  C_\delta \vpran{\Sup_{[0, T]}
  \norm{vf}_{L^1}}
  \vpran{\int_0^T \norm{\vint{v}^{-3/2} \nabla_v \CalM w}^2_{L^2} \dt}^{1/2}
  \vpran{\int_0^T \norm{\CalM w}^2_{L^2} \dt}^{1/2}
\\
& \hspace{-1cm}\leq
  \epsilon_0 \int_0^T \norm{\vint{v}^{-3/2} \nabla_v \CalM w}^2_{L^2} \dt
  + C_\delta \vpran{\Sup_{[0, T]}
  \norm{\vint{v}^3 f}_{L^{3/2}}}
  \int_0^T \norm{\CalM w}^2_{L^2} \dt.
\end{align*}
To bound $\abs{\mathrm{Int}(I_{5,1}^{(1,2)})}$, we apply part $(c)$ in Lemma~\ref{lem:nabla-A-a} and obtain 
\begin{align} \label{bound:A}
   \abs{\mathrm{Int}(I_{5,1}^{(1,2)})}
&\leq
 C \int_0^T \int_{\R^3}
   \abs{\vint{v}^{-3/2}\nabla_v \CalM w} 
  \abs{\int_{\R^3} \frac{|\vint{v_1}^{5/2} \widetilde{v_1 f}|}{|v-v_1|^2} \dv_1} \abs{\CalM w_0} \dv\dt \nn
\\
& \quad \,
  + C \vpran{\Sup_{[0, T]}
  \norm{vf}_{L^1}}
  \int_0^T \norm{\vint{v}^{-3/2} \nabla_v \CalM w}_{L^2} \norm{\CalM w}_{L^2} \dt,
\end{align}
where the first term satisfies a similar bound as the first term on the right-hand side of~\eqref{boudn:I-12-12-interm}. Therefore, 
\begin{align*}
  \abs{\mathrm{Int}(I_{5,1}^{(1,2)})}
\leq
  \epsilon_0 \int_0^T \norm{\vint{v}^{-3/2} \nabla_v \CalM w}^2_{L^2} \dt
  + C_\delta \vpran{\Sup_{[0, T]}
  \norm{\vint{v}^3 f}_{L^{3/2}}}
  \int_0^T \norm{\CalM w}^2_{L^2} \dt,
\end{align*}
which, combined with the estimates for $\mathrm{Int}(I_{5,1}^{(1,1)})$, implies that $\mathrm{Int}(I_{5,1}^{(1)})$, thus $\mathrm{Int}(I_{5,1})$, satisfies a similar bound. Finally, we note that $\mathrm{Int}(I_{5,2})$ has a similar and more regular structure as $\mathrm{Int}(I_{5,2}^{(1)})$. Thus, $\mathrm{Int}(I_5)$ satisfies the desired bound in~\eqref{bound:I-5-6}.
\end{proof}

\section{Uniqueness}

In this section, we apply the {\it a priori} estimates in Section~\ref{Sec:a-priori} to prove the main uniqueness theorem. First, we consider the $L^{3/2}$-solutions from~\cite{GGL25}.
\begin{thm} \label{thm:L3-2-Soln}
Suppose $\vint{v}^{k_0} f^{in} \in L^{3/2}(\R^3)$ with $k_0 \geq 5$. Then the global solution to the Landau-Coulomb equation~\eqref{eq:Landau-Coulomb} constructed in~\cite{GGL25} is unique.  
\end{thm}
\begin{proof}
Recall that the solutions found in \cite{GGL25} are classical $\forall t>0$ and  satisfy 
\begin{align*}
 f \geq 0, 
&\qquad 
 \norm{f(t)}_{L^1(\R^3)} = 1,
\qquad 
 \vint{v}^{k_0} f \in C([0, T]; L^{3/2}(\R^3)),
\\
& \vint{v}^{-\frac{3}{2}}
 \vint{v}^{\frac{3k_0}{4}} \nabla f^{\frac{3}{4}} \in L^2((0, T) \times \R^3),
\end{align*}
and $A[f]$ satisfies~\eqref{cond:A-f}. 

Let  $f, g$ be two such solutions with the same initial data and $w = \vint{v}^2 (f- g)$. Then $w$ satisfies~\eqref{eq:w}and is a classical solution to~\eqref{eq:w} $\forall t>0$. Therefore, the {\it a priori} estimates in Section~\ref{Sec:a-priori} are rigorous. Combining Propositions~\ref{prop:I-1}, \ref{prop:I-2}, and \ref{prop:I-3-I-4}, we obtain the following energy estimate: for any $\epsilon_0 > 0$, there exists $C_{\epsilon_0}$ such that
\begin{align*}
& \quad \,
  \Sup_{[0, T]} \norm{\CalM w}^2_{L^2(\R^3)}
\\
&\leq 
  (\epsilon_0 - c_0)
    \int_0^T \norm{\vint{v}^{-3/2} \nabla_v \CalM w}^2_{L^2} \dt
 + \epsilon_0 \, \Sup_{[0, T]} \norm{\CalM w}^2_{L^2}
 + C_{\epsilon_0}
   \int_0^T 
   \norm{\CalM w}^2_{L^2} \dt
\\
& \leq 
  \vpran{\epsilon_0
    -c_0}
    \int_0^T \norm{\vint{v}^{-3/2} \nabla_v \CalM w}^2_{L^2} \dt
 + \vpran{\epsilon_0 + C_{\epsilon_0} T} \, \Sup_{[0, T]} \norm{\CalM w}^2_{L^2}, 
\end{align*}
where the constant $C_{\epsilon_0}$  depends on $\epsilon_0$ as well as the (uniform) bounds of $\norm{\vint{v}^{k_0} f}_{L^\infty(0, T; L^{3/2})}$ and $\norm{\vint{v}^{-\frac{3}{2}}
 \vint{v}^{\frac{3k_0}{4}} \nabla f^{\frac{3}{4}}}_{L^2((0, T) \times \R^3)}$. Thus, we only need to first take $\epsilon_0 = c_0/2$ and then take $T$ small enough to obtain that $\Sup_{[0, T]} \norm{\norm{\CalM w}}^2_{L^2(\R^3)} = 0$, which gives the uniqueness. 
\end{proof}

\begin{rmk}
If the initial data $f(0, v)$ and $g(0, v)$ are different, then our estimate gives a stability result since it is based on an energy method. 
\end{rmk}

\Ni Now we prove the uniqueness of $H$-solutions, as stated in the Main Theorem in the introduction.

\begin{thm} \label{thm:unique-H}
Suppose $f$ is a nonnegative $H$-solution to~\eqref{eq:Landau-Coulomb} satisfying the extra regularity
\begin{align} \label{cond:f-reg}
  \vint{v}^{k_0} f \in C([0, T]; L^{3/2}(\R^3)),
\qquad  \vint{v}^{-3/2}\vint{v}^{\frac{3k_0}{4}} \nabla f^{\frac{3}{4}} \in L^2((0, T) \times \R^3),
\qquad
  k_0 > 5.
\end{align}
Then $f$ is unique. 
\end{thm}
\begin{proof}
Since it has been shown in the proof of Theorem~\ref{thm:L3-2-Soln} that the energy estimate for $\CalM w$ implies uniqueness, we are left with the justification of the use of the test function $\phi = \vint{v}^2 \CalM^2 w$ in the weak formulation of the $H$-solution. Recall the definition of the $H$-solution: for any test function $\phi \in C_c^2([0, T] \times \R^3)$, it holds that 
\begin{align*} 
& -\int_{\R^3} f^{in} \phi(0, v) \dv
- \int_0^T \int_{\R^3}
   f \del_t \phi \dv \dt
\\
& \hspace{3cm}
  = \int_0^T \int_{\R^3}
   A[f] f : \nabla^2 \phi \dv \dt 
  + \int_0^T \int_{\R^3}
     \nabla a[f] f \cdot \nabla \phi
     \dv \dt.
\end{align*}
Since the bounds of $f$ in~\eqref{cond:f-reg} are global in $v$, we only need $\phi \in C^2([0, T] \times \R^3)$. Suppose $f, g$ are both $H$-solutions with the same initial data. Then their difference $w = \vint{v}^2 (f - g)$ satisfies
\begin{align} \label{eq:H-soln-w}
  - \int_0^T \int_{\R^3}
   w \, \del_t \phi \dv \dt
&  = \int_0^T \int_{\R^3}
   A[f] w : \nabla^2 \phi \dv \dt 
  + \int_0^T \int_{\R^3}
     \nabla a[f] w \cdot \nabla \phi
     \dv \dt  \nn
\\
& \quad \,
  + \int_0^T \int_{\R^3}
   A[w] g : \nabla^2 \phi \dv \dt 
  + \int_0^T \int_{\R^3}
     \nabla a[w] g \cdot \nabla \phi
     \dv \dt, 
\end{align}
and $w$ satisfies the same regularity in~\eqref{cond:f-reg} with a reduced weight $k_0 - 2$. 
Our goal is to show that the energy inequality for $\CalM w$, and hence the {\it a priori} estimates in the previous section, are rigorous.

We first show that, given the regularity
\begin{align} 
  \vint{v}^{3} f \in C([0, T]; L^{3/2}(\R^3)),
\qquad
  \vint{v}^{-3} f \in L^1(0, T; L^3(\R^3)), \label{reg:f-1}
\\
  \vint{v}^{3} g \in C([0, T]; L^{3/2}(\R^3)), 
\qquad
  \vint{v}^{-3} g \in L^1(0, T; L^3(\R^3)), \label{reg:g-1}
\end{align}
every term on the right-hand side of~\eqref{eq:H-soln-w} is well-defined with $\phi = \CalM^2 w$. The second condition in~\eqref{reg:f-1} and~\eqref{reg:g-1} is guaranteed for $H$-solutions, see Theorem 1 in \cite{D14}. First, by HLS,  
\begin{align*}
& \norm{A[f] (t)}_{L^6(\R^3)}
\leq
  C \norm{f(t)}_{L^{6/5} (\R^3)}
\leq 
  C \norm{f(t)}_{L^1 (\R^3)}^{1/2}
  \norm{f}_{L^{3/2} (\R^3)}^{1/2},
\\
 & \norm{\nabla a[f] (t)}_{L^6(\R^3)}
\leq
  C \norm{f (t)}_{L^{2} (\R^3)}
\leq
  C \norm{\vint{v}^{3} f}_{L^{3/2} (\R^3)}^{1/2}
  \norm{\vint{v}^{-3} f}_{L^{3} (\R^3)}^{1/2}.
\end{align*}
Moreover, for $\phi = \CalM^2 w$, we have
\begin{align*}
 & \norm{\nabla^2 \phi (t)}_{L^6(\R^3)}
 = \norm{\nabla^2 \CalM^2 w (t)}_{L^6(\R^3)}
\leq
 C \norm{\CalM w}_{L^6(\R^3)}
\\
& \hspace{6cm}
\leq
 C \norm{\CalM^{3/4} w}_{L^6(\R^3)}
\leq 
 C \norm{w}_{L^{3/2}(\R^3)},
\\
 & \norm{\nabla \phi}_{L^6(\R^3)}
 = \norm{\nabla \CalM^2 w}_{L^6(\R^3)}
\leq 
  C \norm{\CalM w}_{L^2(\R^3)}
\\
& \hspace{6cm}
\leq
 C \norm{\CalM^{1/4} w}_{L^2(\R^3)}  
\leq
 C \norm{w}_{L^{3/2}(\R^3)},
\end{align*}
where we have used the Sobolev embedding in $\R^3$:
\begin{align*}
  W^{\frac{3}{2}, \frac{3}{2}}(\R^3)
\hookrightarrow L^6 (\R^3), 
\qquad
  W^{\frac{1}{2}, \frac{3}{2}}(\R^3)
\hookrightarrow L^2 (\R^3). 
\end{align*}
Therefore,
\begin{align*}
& \quad \,
  \abs{\int_0^T \int_{\R^3}
   A[f] w : \nabla^2 \phi \dv \dt} 
\leq
 \int_0^T \int_{\R^3}
   \abs{A[f]} |w| |\nabla^2 \phi| \dv \dt
\\
&\leq
  \int_0^T
  \norm{A[f]}_{L^6(\R^3)}
  \norm{w}_{L^{3/2}(\R^3)}
  \norm{\nabla^2 \phi}_{L^6(\R^3)} \dt
\\
& \leq
 C T \sup_{[0, T]}
 \vpran{\norm{f(t)}_{L^1 (\R^3)}^{1/2}
  \norm{f}_{L^{3/2} (\R^3)}^{1/2}}
  \sup_{[0, T]}
 \vpran{\norm{w}_{L^{3/2} (\R^3)}^2}
< \infty,
\end{align*}
and
\begin{align*}
& \quad \,
  \abs{\int_0^T \int_{\R^3}
   \nabla a[f] w \cdot \nabla \phi \dv \dt} 
\leq
 \int_0^T \int_{\R^3}
   \abs{\nabla a[f]} |w| |\nabla \phi| \dv \dt
\\
&\leq
  \int_0^T 
  \norm{\nabla a[f]}_{L^6(\R^3)}
  \norm{w}_{L^{3/2}(\R^3)}
  \norm{\nabla \phi}_{L^6(\R^3)} \dt 
\\
& \leq
  C \sup_{[0, T]}
 \vpran{\norm{\vint{v}^3 f(t)}_{L^1 (\R^3)}^{1/2} \norm{w}_{L^{3/2} (\R^3)}^2}
 \int_0^T
   \norm{\vint{v}^{-3} f}_{L^{3} (\R^3)}^{1/2} \dt
\\
& \leq
 C \sqrt{T}
 \sup_{[0, T]}
 \vpran{\norm{\vint{v}^3 f(t)}_{L^1 (\R^3)}^{1/2} \norm{w}_{L^{3/2} (\R^3)}^2}
 \vpran{\int_0^T
   \norm{\vint{v}^{-3} f}_{L^{3} (\R^3)} \dt}^{1/2}
< \infty,
\end{align*}
Similar arguments apply to the third and fourth terms of (\ref{eq:H-soln-w}), since $w$, $f$ and $g$ enjoy the same regularity. 

Once we have (\ref{eq:H-soln-w}), the rigorous argument can then be performed by density argument via mollifiers. In particular, one considers the sequence of test functions
\begin{align*}
  \phi_{\delta, \epsilon}
= \mathcal{N}_\delta^2 \CalM^2 \CalM_\epsilon^3 w,
\end{align*}
where $\mathcal{N}_{\delta}$ is a convolution mollifier in time and $\CalM_\epsilon
= (I - \epsilon \Delta_v)^{-1}$. The energy estimate with $\phi = \CalM^2 w$ is then justified by passing $\delta, \epsilon \to 0$ via the Lebesgue Dominated Convergence theorem.
\end{proof}

We remark that the condition $\langle v \rangle^{m} \nabla f^{\frac{3}{4}}$  is not necessary to justify (\ref{eq:H-soln-w}) with $\phi = \CalM^2 w$. It is, however, necessary to show uniqueness. In this regard, we close with a final remark:  the {\it a priori} estimate for $\nabla f^{\frac{3}{4}}$ can be derived from the bound of $f$ in $C([0, T]; L^{3/2}_{k}(\R^3))$ for $k > \frac{18}{5}$. Nevertheless, as is common for nonlinear PDEs, such {\it a priori} estimates do not automatically apply to weak solutions.  Therefore, in Theorem~\ref{thm:unique-H} we must include the condition on $\nabla f^{\frac{3}{4}}$ . 

\begin{prop}
Suppose $f$ is a sufficiently smooth nonnegative solution satisfying 
\begin{align*}
  \vint{v}^{k} f \in C([0, T]; L^{3/2}(\R^3)),
\qquad
  k > \frac{18}{5}.
\end{align*}
Then $\vint{v}^{-3/2}\vint{v}^{\frac{3k}{4}} \nabla f^{\frac{3}{4}} \in L^2((0, T) \times \R^3)$. 
\end{prop}

\begin{proof}
Our goal is to derive the bound of $\vint{v}^{k-3} \nabla f^{\frac{3}{4}}$ in terms of $\norm{\vint{v}^{k} f}_{L^\infty(0, T; L^{3/2})}$. Perform the $L^{3/2}$-energy estimate for $\vint{v}^{k}$ by using $\vint{v}^{\frac{3k}{2}} \sqrt{f}$ as the test function for~\eqref{eq:Landau-Coulomb}. We have
\begin{align} \label{bound:f-energy}
  \frac{\rm d}{\dt}
\norm{\vint{v}^k f}^{\frac{3}{2}}_{L^{\frac{3}{2}}}
& = - \int_{\R^3}
    \nabla \vpran{\vint{v}^{\frac{3k}{2}} \sqrt{f}} \cdot A[f] \cdot \nabla f \dv
    - \int_{\R^3}
    \nabla \vpran{\vint{v}^{\frac{3k}{2}} \sqrt{f}} \cdot \nabla a[f] \, f \dv  \nn
\\
& \hspace{-1cm}
   = - \frac{8}{9} \int_{\R^3}
    \vpran{\vint{v}^{\frac{3k}{4}} \nabla f^{\frac{3}{4}}} \cdot A[f] \cdot \vpran{\vint{v}^{\frac{3k}{4}} \nabla f^{\frac{3}{4}}} \dv
   + \frac{1}{3}
    \int_{\R^3}
     \vint{v}^{\frac{3k}{2}} f^{\frac{5}{2}} \dv  \nn
\\
& \hspace{-0.5cm} 
  - \frac{4}{3} \int_{\R^3}
    \nabla \vpran{\vint{v}^{\frac{3k}{2}}} f^{\frac{3}{4}} \cdot A[f] \cdot \nabla f^{\frac{3}{4}} \dv
  - \int_{\R^3}
   \nabla \vpran{\vint{v}^{\frac{3k}{2}}} f^{\frac{3}{2}} \cdot \nabla a[f] \dv \nn
\\
& \hspace{-1cm} 
  =: E_1 + E_2 + E_3 + E_4. 
\end{align}
where to derive $E_2$, we have used $-\Delta a[f] = f$. Among these terms, $E_1$ and $E_2$ are the leading orders and $E_3, E_4$ are of lower order. 
The first term $E_1$ provides the diffusion:
\begin{align*}
  E_1 
\leq 
  - c_0 \int_{\R^3}
    \vint{v}^{-3} \abs{\vint{v}^{\frac{3k}{4}} \nabla f^{\frac{3}{4}}}^2 \dv_1.
\end{align*}
Separating $f = f_\delta + \tilde f$ in $E_2$. Then we have
\begin{align*}
  E_2
\leq 
 C  \int_{\R^3}
     \vint{v}^{\frac{3k}{2}} f^{\frac{5}{2}}_\delta \dv
 + C \int_{\R^3}
     \vint{v}^{\frac{3k}{2}} \tilde f^{\frac{5}{2}} \dv
=: E_{2,1} + E_{2,2}.
\end{align*}
For any $\delta > 0$, the first term $E_{2,1}$ can be bounded directly as
\begin{align*}
  E_{2,1}
\leq
 C \norm{f_\delta}_{L^\infty}
 \norm{\vint{v}^k f_\delta}_{L^{\frac{3}{2}}}^{\frac{3}{2}}
\leq
 \frac{C}{\delta} \norm{f}_{L^1}
 \norm{\vint{v}^k f}_{L^{\frac{3}{2}}}^{\frac{3}{2}}
\leq 
  C \norm{\vint{v}^k f}_{L^{\frac{3}{2}}}^{\frac{3}{2}}.
\end{align*}
Apply H\"{o}lder's inequality and Sobolev embedding to $E_{2,2}$ and we get
\begin{align*}
  E_{2,2}
&\leq 
 C \norm{\vint{v}^{-\frac{3}{2}} \vint{v}^{\frac{3k}{4}} f}_{L^{\frac{9}{2}}}^{\frac{3}{2}}
 \norm{\vint{v}^{\frac{3k}{8} + \frac{9}{4}} \tilde f}_{L^{\frac{3}{2}}}
\\
&\leq
 C \vpran{\norm{\vint{v}^{-3/2} \vint{v}^{\frac{3k}{4}} \nabla f^{\frac{3}{4}}}_{L^2}^{2}
 + \norm{\vint{v}^k f}_{L^{3/2}}^{3/2}}
 \norm{\vint{v}^{k} \tilde f}_{L^{\frac{3}{2}}}
\\
& \leq 
 \frac{c_0}{2}
 \norm{\vint{v}^{-3/2} \vint{v}^{\frac{3k}{4}} \nabla f^{\frac{3}{4}}}_{L^2}^{2}
 + C \norm{\vint{v}^k f}_{L^{\frac{3}{2}}}^{\frac{5}{2}},
\end{align*}
where we have chosen $\delta$ small enough such that 
\begin{align*}
  C \norm{\vint{v}^{k} \tilde f}_{L^{\frac{3}{2}}}
\leq 
  \frac{c_0}{2},
\end{align*}
and $k$ large enough such that $\frac{3k}{8} + \frac{9}{4} \leq k$ or equivalently, $k \geq \frac{18}{5}$. 
Summarizing these estimates, we have
\begin{align*}
  E_{2}
\leq
 C \norm{\vint{v}^k f}_{L^{\frac{3}{2}}}^{\frac{3}{2}}
 + C \norm{\vint{v}^k f}_{L^{\frac{3}{2}}}^{\frac{5}{2}}
 + \frac{c_0}{2} \norm{\vint{v}^{-3} \nabla f^{\frac{3}{4}}}_{L^2}^{2}
\leq
 C_f + \frac{c_0}{2} \norm{\vint{v}^{-3} \nabla f^{\frac{3}{4}}}_{L^2}^{2},
\end{align*}
where $C_f$ depends on $\sup_{[0, T]} \norm{\vint{v}^k f}_{L^{\frac{3}{2}}}$.

Now we show the bounds for the lower order terms $E_3$ and $E_4$. First, by integration by parts,  
\begin{align*}
  E_3
&= - \frac{4}{3} \int_{\R^3}
    \nabla \vpran{\vint{v}^{\frac{3k}{2}}} f^{\frac{3}{4}} \cdot A[f] \cdot \nabla f^{\frac{3}{4}} \dv
\\
&= \frac{2}{3} \int_{\R^3}
  \nabla^2 \vpran{\vint{v}^{\frac{3k}{2}}} : A[f] \, f^{\frac{3}{2}} \dv
  + \frac{2}{3} \int_{\R^3}
  \nabla \vpran{\vint{v}^{\frac{3k}{2}}} \cdot (\nabla \cdot A[f]) \, f^{\frac{3}{2}} \dv
=: E_{3,1} + E_{3,2},
\end{align*}
where, by Cauchy-Schwarz and HLS, 
\begin{align*}
  |E_{3,1}|
\leq
  C \norm{A[f]}_{L^6}
  \norm{\vint{v}^{\frac{3k}{2}-2} f}_{L^{\frac{9}{5}}}^{\frac{3}{2}}
&\leq
 C \norm{f}_{L^{\frac{6}{5}}}
 \norm{\vint{v}^k f}_{L^{\frac{3}{2}}}^{\frac{9}{8}}
 \norm{\vint{v}^{-3/2} \vint{v}^k f}_{L^{\frac{9}{2}}}^{\frac{3}{8}}
\\
& \leq 
 C \norm{\vint{v}^k f}_{L^{\frac{3}{2}}}^{\frac{17}{8}}
 \norm{\vint{v}^{-3/2} \vint{v}^k f}_{L^{\frac{9}{2}}}^{\frac{3}{8}},
\end{align*}
for $k > 1/2$. By Young's inequality, for any $\epsilon_0 > 0$,
\begin{align*}
  |E_{3,1}|
&\leq
 C_{\epsilon_0} \norm{\vint{v}^k f}_{L^{\frac{3}{2}}}^{\frac{17}{6}}
 + \epsilon_0 \norm{\vint{v}^{-3/2} \vint{v}^k f}_{L^{\frac{9}{2}}}^{\frac{3}{2}}
\\
&\leq 
 C_f
 + \epsilon_0 \norm{\vint{v}^{-3/2} \vint{v}^k \nabla f^{\frac{3}{4}}}^2_{L^2},
\end{align*}
where again, $C_f$ depends on $\sup_{[0, T]} \norm{\vint{v}^k f}_{L^{\frac{3}{2}}}$.

Finally, By Lemma~\ref{lem:nabla-A-a} part (a) and Young's inequality,
\begin{align*}
  |E_{3,2}|
&\leq
 C \norm{f}_{L^1}
 \norm{\vint{v}^k f}_{L^{\frac{3}{2}}}^{\frac{3}{2}}
 + C \norm{\int_{\R^3} 
 \frac{|\vint{v_1}^{2} f (v_1)|}{|v-v_1|^2} \dv_1}_{L^3} 
 \norm{\vint{v}^{\frac{3k}{2}-3} f^{\frac{3}{2}}}_{L^{\frac{3}{2}}}
\\
& \leq
  C \norm{\vint{v}^k f}_{L^{\frac{3}{2}}}^{\frac{3}{2}}
  + C \norm{\vint{v}^2 f}_{L^{\frac{3}{2}}}
  \norm{\vint{v}^{k} f}_{L^{\frac{3}{2}}}^{\frac{3}{4}}
  \norm{\vint{v}^{-\frac{3}{2}} \vint{v}^k f}_{L^{\frac{9}{2}}}^{\frac{3}{4}}
\\
&\leq
 C_f
 + \epsilon_0 \norm{\vint{v}^{-\frac{3}{2}} \vint{v}^k \nabla f^{\frac{3}{4}}}_{L^2}^2, 
\end{align*}
with $C_f$ depending on $\sup_{[0, T]} \norm{\vint{v}^k f}_{L^{\frac{3}{2}}}$. 
Therefore, 
\begin{align*}
  |E_3|
\leq 
 2\epsilon_0 \norm{\vint{v}^{-\frac{3}{2}} \vint{v}^k \nabla f^{\frac{3}{4}}}_{L^2}^2
 + C_{\epsilon_0, f},
\end{align*}
where $C_{f, \epsilon_0}$ depends on $\epsilon$ and $\sup_{[0, T]} \norm{\vint{v}^k f}_{L^{\frac{3}{2}}}$. The last term $E_4$ satisfies that $E_4 = \frac{1}{2} E_{3,2}$. Therefore, it satisfies a similar bound as $E_3$. Applying the estimates for $E_1, \cdots, E_4$ to~\eqref{bound:f-energy} and integrating in $t$, we have
\begin{align*}
  \norm{\vint{v}^k f}^{\frac{3}{2}}_{L^{\frac{3}{2}}}
- \norm{\vint{v}^k f^{in}}^{\frac{3}{2}}_{L^{\frac{3}{2}}}
\leq 
 -\frac{c_0}{4} \int_0^t \norm{\vint{v}^{-3} \nabla f^{\frac{3}{4}}}_{L^2}^{2} \ds
+ t C_{f} , 
\end{align*}
by choosing $\epsilon_0 = c_0/8$. This gives
\begin{align*}
 \int_0^T \norm{\vint{v}^{-3} \vint{v}^{\frac{3k}{4}}\nabla f^{\frac{3}{4}}}_{L^2}^{2} \ds
&\leq 
  C_f T,
< \infty,
\end{align*}
where $C_f$ depends on $\sup_{[0, T]} \norm{\vint{v}^k f}_{L^{\frac{3}{2}}}$ algebraically.
\end{proof}

\appendix

\section{Proof of Lemma~\ref{ineq:basic}}

In this appendix, we prove parts $(e)$ and $(f)$ in Lemma~\ref{ineq:basic}. 

\begin{proof}
$(e)$ Denote $M (\xi)$ as the symbol of $\CalM$. Then, for any $\beta > 0$, $M^\beta$ satisfies
\begin{align*}
  \abs{\del_\xi^{\alpha} M^\beta(\xi)}
\leq
  C_\alpha \vint{\xi}^{-|\alpha| - 2\beta}
\leq
  C_\alpha \vint{\xi}^{-|\alpha|},
\end{align*}
for all multi-index $\alpha$. Therefore, Mikhlin's multiplier theorem implies that $\CalM^\beta$ is bounded on $L^p$ for any $1 < p < \infty$. Similarly, the symbols of $\nabla \CalM$ and $\nabla^2 \CalM$ satisfy
\begin{align*}
  \abs{\del_\xi^\alpha (\xi M(\xi))}
\leq 
  C_\alpha \vint{\xi}^{-|\alpha|-1}
\leq
  C_\alpha \vint{\xi}^{-|\alpha|},
\qquad
  \abs{\del_\xi^\alpha (\xi \otimes \xi \, M(\xi))}
\leq 
  C_\alpha \vint{\xi}^{-|\alpha|}.
\end{align*}
Therefore, $\nabla \CalM$, $\nabla^2 \CalM$ are bounded on $L^p$ for any $1 < p < \infty$. Thus, the last two inequalities in $(e)$ hold. For the first inequality, note that for $f \in L^p$, we have $\CalM f \in W^{2,p}$ by the third inequality. Therefore, by Sobolev embedding, if $p < \frac{n}{2}$, then $\dot{W}^{2, p} \hookrightarrow L^q$ with $\frac{1}{q} = \frac{1}{p} - \frac{2}{n}$. The full range of $q$ follows by interpolation, since the second inequality gives $\CalM f \in L^p$. 

\smallskip

\Ni $(f)$ The inequalities in part $(f)$ can be shown by commutator estimates for pseudo-differential operators. However, given the explicit integral representation of 
$\CalM$, we give direct proofs. 

To prove~\eqref{lem-ineq:f-1}, denote the commutator as
\begin{align*}
  \CalT f
:= \vint{v}^\alpha \CalM f
  - \CalM (\vint{v}^\alpha f). 
\end{align*}
Our goal is to prove
\begin{align} \label{est:T-1}
  \norm{\CalT f}_{L^p}
\leq 
  C \norm{\vint{v}^\alpha \CalM f}_{L^p},
\end{align}
and
\begin{align} \label{est:T-2}
  \norm{\CalT f}_{L^p(\R^n)}
\leq 
  C \norm{\CalM(\vint{v}^\alpha f)}_{L^p(\R^n)}.
\end{align}
To prove~\eqref{est:T-1}, we rewrite $\CalT f$ as
\begin{align*}
  \CalT f
&= \CalM (I - \Delta) (\vint{v}^\alpha \CalM f)
  - \CalM (\vint{v}^\alpha f)
= -\CalM
  \vpran{(\Delta \vint{v}^\alpha) \CalM f}
  - 2 \CalM \vpran{\nabla \vint{v}^\alpha \cdot \nabla \CalM f}
\\
&= -3 \CalM
  \vpran{(\Delta \vint{v}^\alpha) \CalM f}
  -2 \nabla \CalM \cdot
    ((\nabla \vint{v}^\alpha) \, \CalM f)
\end{align*}
where 
\begin{align*}
  \norm{\CalM
  \vpran{(\Delta \vint{v}^\alpha) \CalM f}}_{L^p}
\leq
 C_p 
 \norm{(\Delta \vint{v}^\alpha) \, \CalM f}_{L^p}
\leq
 C_p \norm{\vint{v}^{\alpha-2} \CalM f}_{L^p}
\leq
 C_p \norm{\vint{v}^{\alpha} \CalM f}_{L^p},
\end{align*}
and by part $(e)$,
\begin{align*}
 \norm{\nabla\CalM \cdot \vpran{(\nabla \vint{v}^\alpha) \,\, \CalM f}}_{L^p}
\leq
 C_p \norm{(\nabla \vint{v}^\alpha) \,\, \CalM f}_{L^p}
\leq
 C_{p,\alpha} \norm{\vint{v}^\alpha \CalM f}_{L^p}.
\end{align*}
Therefore, \eqref{est:T-1} holds, which implies that
\begin{align*}
 \norm{\CalM \vpran{\vint{v}^{\alpha} f}}_{L^p(\R^n)}
\leq
  C \norm{\vint{v}^{\alpha} \CalM f}_{L^p(\R^n)}.
\end{align*}
To prove~\eqref{est:T-2}, we rewrite $\CalT f$ as
\begin{align*}
 \CalT f
&= \vint{v}^\alpha \CalM \vpran{\vint{v}^{-\alpha} (I - \Delta) \CalM (\vint{v}^\alpha f)}
  - \CalM (\vint{v}^\alpha f)
\\
&= \vint{v}^\alpha \CalM \vpran{(\Delta \vint{v}^{-\alpha}) \CalM (\vint{v}^\alpha f}
 + 2\vint{v}^\alpha \CalM \vpran{
   \nabla \vint{v}^{-\alpha} \cdot \nabla (\CalM (\vint{v}^\alpha f)}
\\
&= -\vint{v}^\alpha \CalM \vpran{(\Delta \vint{v}^{-\alpha}) \CalM (\vint{v}^\alpha f}
 + \vint{v}^\alpha \nabla\CalM \cdot \vpran{
   (\nabla \vint{v}^{-\alpha})  (\CalM (\vint{v}^\alpha f)}.
\end{align*}
Therefore, if we define
\begin{align*}
  \CalT_1 G
= \vint{v}^\alpha \CalM \vpran{(\Delta \vint{v}^{-\alpha}) G},
\qquad
  \CalT_2 G
= \vint{v}^\alpha \nabla\CalM \cdot \vpran{
   (\nabla \vint{v}^{-\alpha}) G},
\end{align*}
then it suffices to prove that $\CalT_1, \CalT_2$ are both bounded on $L^p$. By the integral form of $\CalM$, we have
\begin{align*}
  \CalT_1 G
&= C \vint{v}^\alpha 
  \int_{\R^3}
  \frac{e^{-|v-v_1|}}{|v-v_1|}
  \vpran{\Delta \vint{v_1}^{-\alpha}} \, G(v_1) \dv_1
\\
& \leq
  C \vint{v}^\alpha 
  \int_{\R^3}
  \frac{e^{-|v-v_1|}}{|v-v_1|}
  \vint{v_1}^{-\alpha-2} \, |G(v_1)| \dv_1
\\
& = C \vint{v}^\alpha 
  \int_{|v_1| \leq \frac{1}{2}\vint{v}}
  + \, C \vint{v}^\alpha 
  \int_{|v_1| \geq 2\vint{v}}
  + \, C \vint{v}^\alpha 
  \int_{\frac{1}{2}\vint{v} < |v_1| < 2 \vint{v} },
\end{align*}
where over the three integration domains, we have
\begin{align*}
 &e^{-\frac{1}{2}|v-v_1|}
 \vint{v}^\alpha \vint{v_1}^{-\alpha -2}
\leq
 e^{-\frac{1}{4}|v|} \vint{v}^{|\alpha|} \vint{v_1}^{|\alpha|}
< C, 
\qquad
 |v_1| \leq \tfrac{1}{2} \vint{v},
\\
 & e^{-\frac{1}{2}|v-v_1|}
 \vint{v}^\alpha \vint{v_1}^{-\alpha -2}
\leq
 e^{-\frac{1}{4}|v_1|} \vint{v}^{|\alpha|} \vint{v_1}^{|\alpha|}
< C, 
\qquad
 |v_1| \geq 2 \vint{v},
\\
 & \vint{v}^\alpha \vint{v_1}^{-\alpha -2}
\leq
 C,
\hspace{6cm}
  \tfrac{1}{2}\vint{v} < |v_1| < 2 \vint{v}.
\end{align*}
Therefore, over all the subdomains, we have
\begin{align*}
 \abs{\CalT_1 G}
\leq
 C \int_{\R^3}
  \frac{e^{-\frac{1}{2}|v-v_1|}}{|v-v_1|}
  |G(v_1)| \dv_1.
\end{align*}
Since the kernel $\frac{e^{-\frac{1}{2} |v_1|}}{|v_1|} \in L^1(\R^3)$, by Young's inequality, we have 
\begin{align*}
  \norm{\CalT_1 G}_{L^p(\R^3)} \leq C \norm{G}_{L^p(\R^3)},
\qquad
  \forall \, 1 < p < \infty.
\end{align*}
The $L^p \to L^p$ bound of $\CalT_2$ is shown in a similar way:
\begin{align*}
  \abs{\CalT_2 G}
&\leq
  C \vint{v}^\alpha 
  \nabla_v \int_{\R^3}
  \frac{e^{-|v-v_1|}}{|v-v_1|}
  \cdot 
  \vpran{\nabla \vint{v_1}^{-\alpha}} \, G(v_1) \dv_1
\\
& \leq
  C \vint{v}^\alpha 
  \int_{\R^3}
  \frac{e^{-|v-v_1|}}{|v-v_1|^2}
  \vint{v_1}^{-\alpha-1} \, |G(v_1)| \dv_1.
\end{align*}
By a similar estimate as for $\CalT_1$, we have
\begin{align*}
  e^{-\frac{1}{2}|v-v_1|}
 \vint{v}^\alpha \vint{v_1}^{-\alpha -2}
< C,
\qquad 
 \forall \, v, v_1 \in \R^3.
\end{align*}
Moreover, the kernel $\frac{e^{-\frac{1}{2} |v_1|}}{|v_1|^2} \in L^1(\R^3)$. Therefore, $\CalT_2$ is also a bounded operator on $L^p(\R^3)$ for any $1 < p < \infty$, which shows
\begin{align*}
  \norm{\vint{v}^\alpha \CalM f}_{L^p}
\leq
  \norm{\CalM (\vint{v}^\alpha f)}_{L^p} + \norm{\CalT f}_{L^p}
\leq
  C \norm{\CalM (\vint{v}^\alpha f)}_{L^p}.
\end{align*}

Next, we show details for proving~\eqref{lem-ineq:f-2} with $|\alpha|=2$. The case of $|\alpha| = 1$ can be shown in a similar way. To this end, we prove that
\begin{align*}
  \norm{\vint{v}^\beta \nabla^2\CalM f}_{L^p(\R^3)}
\leq
  C \norm{\vint{v}^\beta f}_{L^p(\R^3)},
\end{align*}
for any $\beta \in \R $ and $1 < p < \infty$. 
Rewrite the left-hand side as
\begin{align*}
  \norm{\vint{v}^\beta \nabla^2\CalM \vpran{\vint{v}^{-\beta} \vpran{\vint{v}^\beta f}}}_{L^p(\R^3)}. 
\end{align*}
Then our goal becomes showing that the operator $\vint{v}^\beta \nabla^2\CalM \vpran{\vint{v}^{-\beta} \cdot}$ is bounded on $L^p(\R^3)$ for any $1 < p < \infty$. For any $h$ sufficiently smooth, by part $(e)$, we have
\begin{align} \label{A-ineq:1}
  \norm{\vint{v}^\beta \nabla^2\CalM \vpran{\vint{v}^{-\beta} h}}_{L^p(\R^3)}
&= \norm{\vint{v}^\beta \CalM \vpran{\nabla^2 \vpran{\vint{v}^{-\beta} h}}}_{L^p(\R^3)}  \nn
\\
& \leq
 C \norm{\CalM \vpran{\vint{v}^\beta \nabla^2 \vpran{\vint{v}^{-\beta} h}}}_{L^p(\R^3)},
\end{align}
where 
\begin{align*}
 \vint{v}^\beta \nabla^2 \vpran{\vint{v}^{-\beta} h}
&= \nabla^2 h
  - 2 \nabla \vint{v}^\beta
    \otimes \nabla \vpran{\vint{v}^{-\beta} h}
  - \vpran{\nabla^2 \vint{v}^\beta}
  \vint{v}^{-\beta} h
\\
&= \nabla^2 h
   - 2 \nabla \vpran{(\nabla\vint{v}^\beta)\vint{v}^{-\beta} h}
   + \vpran{\nabla^2 \vint{v}^\beta}
  \vint{v}^{-\beta} h.
\end{align*}
Therefore, continuing the estimates in~\eqref{A-ineq:1}, we have
\begin{align*}
    \norm{\CalM \vpran{\vint{v}^\beta \nabla^2 \vpran{\vint{v}^{-\beta} h}}}_{L^p(\R^3)}
& \leq
 \norm{\CalM \nabla^2 h}_{L^p(\R^3)}
 + \norm{\CalM \vpran{\vpran{\nabla^2 \vint{v}^\beta}
  \vint{v}^{-\beta} h}}_{L^p(\R^3)}
\\
& \quad \,
  + 2 \norm{\CalM \nabla \vpran{(\nabla\vint{v}^\beta)\vint{v}^{-\beta} h}}_{L^p(\R^3)}
\\
& \hspace{-3.5cm} \leq
  C \norm{h}_{L^p(\R^3)}
  + C \norm{\vpran{\nabla^2 \vint{v}^\beta}
  \vint{v}^{-\beta} h}_{L^p(\R^3)}
 + C \norm{ (\nabla\vint{v}^\beta)\vint{v}^{-\beta} h}_{L^p(\R^3)}
\\
& \hspace{-3.5cm} \leq 
  C \norm{h}_{L^p(\R^3)},
\end{align*}
where the second last step follows from part $(e)$. Hence, by letting $h = \vint{v}^\beta f$, we obtain that
\begin{align*}
  \norm{\vint{v}^\beta \nabla^2\CalM f}_{L^p(\R^3)}
&= \norm{\vint{v}^\beta \nabla^2\CalM \vpran{\vint{v}^{-\beta} \vpran{\vint{v}^\beta f}}}_{L^p(\R^3)}
\leq
  C \norm{\vint{v}^\beta f}_{L^p(\R^3)}. \qedhere
\end{align*}
\end{proof}

\medskip

\subsection*{Acknowledgments} MPG is partially supported by the DMS-NSF 2511625. WS is partially supported by NSERC Discovery Grant R832717. This work was initiated and largely completed while both authors were attending the semester-long program \textit{Kinetic Theory: Novel Statistical, Stochastic and Analytical Methods} at the Simons Laufer Mathematical Sciences Institute (SLMath) in Berkeley, California, in 2025. We greatly appreciate the hospitality and stimulating environment provided by 
SLMath. 


\bibliographystyle{amsxport}
\bibliography{kinetic}

\end{document}